\documentclass{siamltex}

\usepackage{amsfonts} 
\usepackage{amssymb,amsmath} 
\usepackage{graphicx,epstopdf}

\DeclareMathAlphabet{\itbf}{OML}{cmm}{b}{it}
 
\def\EE{\mathbb{E}}
\def\PP{\mathbb{P}}

\def\RR{\mathbb{R}}
\def\NN{\mathbb{N}}
\def\eps{\varepsilon}

\def\by{{{\itbf y}}}
\def\bx{{{\itbf x}}}

\def\Nr{N_{\rm r}}
\def\Ns{N_{\rm s}}
\def\Nq{N_{\rm q}}

\newcommand{\bk}{{\boldsymbol{\kappa}}}
\newcommand{\bxi}{{\boldsymbol{\xi}}}

\newcommand{\ea}{\end{eqnarray}}  
\newcommand{\ba}{\begin{eqnarray}}  
\newcommand{\ee}{\end{equation}}  
\newcommand{\be}{\begin{equation}}  
\newcommand{\ean}{\end{eqnarray*}}  
\newcommand{\ban}{\begin{eqnarray*}}

\begin{document}

\title{Role of scattering in virtual source array imaging}

\author{Josselin Garnier\footnotemark[1] \,   and George Papanicolaou\footnotemark[2]}   

\maketitle

\renewcommand{\thefootnote}{\fnsymbol{footnote}}

\footnotetext[1]{Laboratoire de Probabilit\'es et Mod\`eles Al\'eatoires
\& Laboratoire Jacques-Louis Lions,
Universit{\'e} Paris Diderot,
75205 Paris Cedex 13,
France
{\tt garnier@math.univ-paris-diderot.fr}}
\footnotetext[2]{Mathematics Department,
Stanford University,
Stanford, CA 94305
{\tt papanicolaou@stanford.edu}}

\renewcommand{\thefootnote}{\arabic{footnote}}

\begin{abstract}
We consider imaging in a scattering medium
where the illumination goes through this medium but
there is also an auxiliary, passive receiver array that is 
near the object to be imaged. Instead
of imaging with the source-receiver array on the far side of the
object we image with the data of the passive array
on the near side of the object. The imaging is done with
travel time migration using the cross correlations of the passive
array data. We showed in [J. Garnier and G. Papanicolaou, 
Inverse Problems {28} (2012), 075002] 
that if (i) the source array is infinite, 
(ii) the scattering medium is modeled
by either an isotropic random medium in the paraxial
regime or a randomly layered medium, and (iii) the
medium between the auxiliary array and the object to
be imaged is homogeneous, then imaging with cross
correlations completely eliminates the effects of the random
medium. It is as if we imaged with an active array, instead
of a passive one, near the object. The purpose of this paper
is to analyze the resolution of the image when both the source array
and the passive receiver array are finite. We show with a detailed
analysis that for isotropic
random media in the paraxial regime, imaging not only is not
affected by the inhomogeneities but the resolution can in fact be 
enhanced. This is because the random medium can increase the
diversity of the illumination. We also show analytically that this will not happen
in a randomly layered medium, and there may be some loss of resolution
in this case.
\end{abstract}



\pagestyle{myheadings}
\thispagestyle{plain}
\markboth{J. Garnier and G. Papanicolaou}
{Role of scattering in virtual source array imaging}

\section{Introduction}
\begin{figure}
\begin{center}
\begin{tabular}{c}
\includegraphics[width=5.5cm]{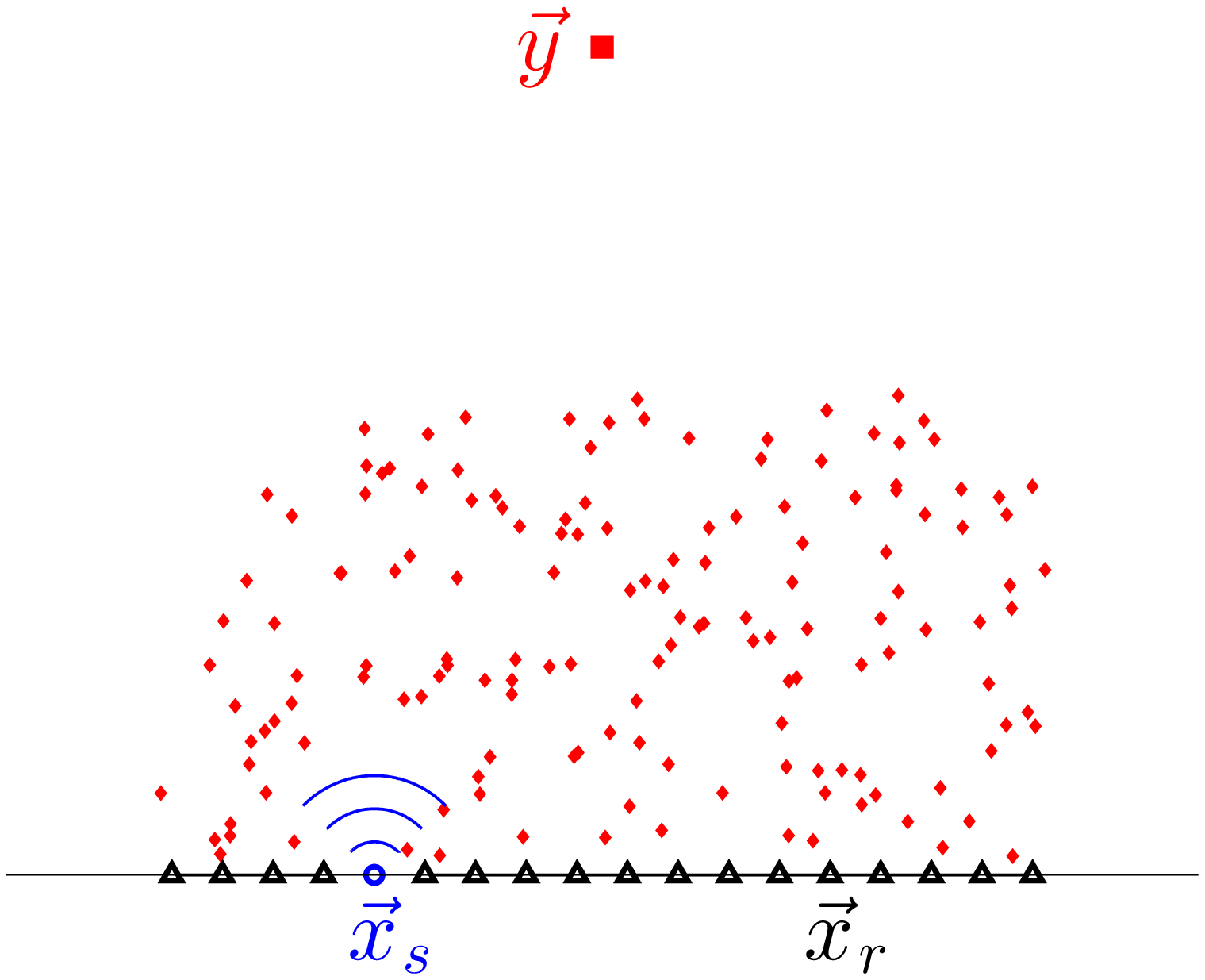} 
\includegraphics[width=5.5cm]{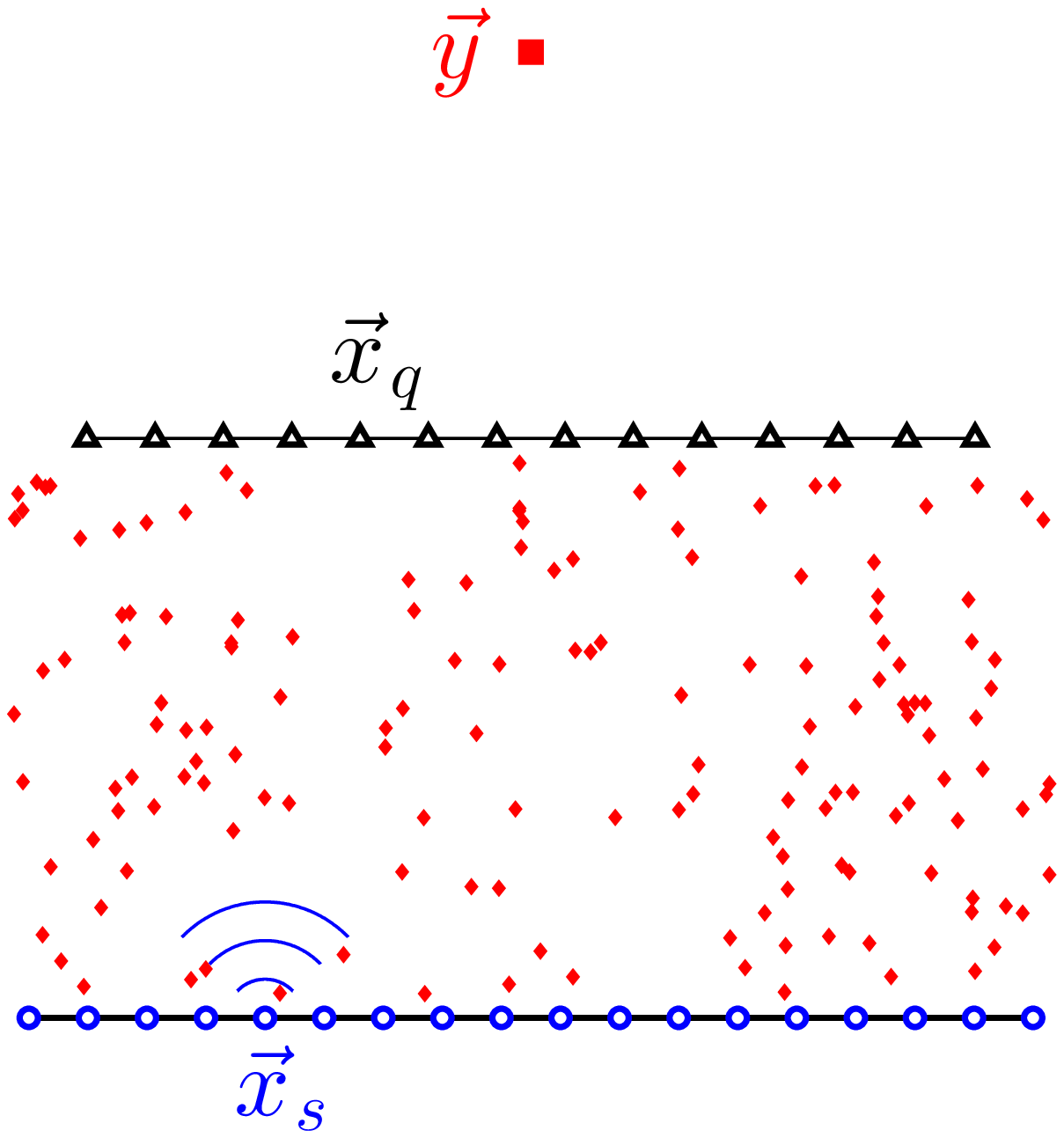} 
\end{tabular}
\end{center}
\caption{Sensor array imaging of a reflector located at $\vec\by$ through a scattering medium. 
Left: Conventional imaging, 
$\vec\bx_s$ is a source,
$\vec\bx_r$ is a receiver.
Right: Use of an auxiliary passive array for imaging through a scattering medium,
$\vec\bx_s$ is a source, $\vec\bx_q$ is a receiver located beyond the scattering medium.}
\label{fig:intro0}
\end{figure}

In conventional active array imaging (see Figure \ref{fig:intro0}, left), 
the sources are on an array located at
$( \vec\bx_s )_{s=1}^{\Ns}$  and the receivers at $( \vec\bx_r )_{r=1}^{\Nr}$, 
where the two arrays are coincident. 
The array response matrix is given by
\begin{equation}
\big\{ p(t,\vec\bx_r;\vec\bx_s), \, t\in \RR, \, r=1,\ldots,\Nr, \, s=1,\ldots,\Ns \big\}
\end{equation}
and consists of the signals recorded by the $r$th receiver when the $s$th source 
emits a short pulse. An image is formed by migrating the array response matrix.
The Kirchhoff migration imaging function \cite{biondi,bleistein} 
at a search point $\vec\by^S$ in the image domain is given by
\begin{equation}
{\cal I}(\vec\by^S) =\frac{1}{\Ns \Nr}  \sum_{r=1}^{\Nr}  \sum_{s=1}^{\Ns} p\Big(
  \frac{ | \vec\bx_s -  \vec\by^S | + | \vec\by^S -\vec\bx_{r} |}{c_0} ,\vec\bx_r; \vec\bx_s \Big) .
\end{equation}
Here $|\vec\bx-\vec\by|/c_0$ is a computed travel time between the points $\vec\bx$ and $\vec\by$,
corresponding to a homogeneous medium with 
constant propagation speed $c_0$. In this case the images produced by Kirchhoff migration, that is,
the peaks of the imaging function ${\cal I}(\vec\by^S) $,
 can be analyzed easily. For a point reflector the range resolution is $c_0/B$, where
$B$ is the bandwidth of the probing pulse, and the cross range resolution is $\lambda_0 L/a$, where $\lambda_0$ is the central wavelength of the pulse, $L$ is the distance from the array to the reflector, and $a$ is the size
of the array. These are the well-known Rayleigh resolution limits \cite{elmore}.
When, however, the medium is inhomogeneous then migration may not work well. 
In weakly scattering media the images can be stabilized statistically by using
coherent interferometry \cite{borcea05,borcea06a,borcea06b,borcea11a,borcea11b}, which is a special correlation-based imaging method.
Statistical stability here means high signal-to-noise ratio for the imaging function.
In strongly scattering media we may be able to obtain an image by using special
signal processing methods \cite{borcea09} but often we cannot get any image at all because the
coherent signals from the reflector received at the array are very weak compared to the
backscatter from the medium. 

Let as also consider the possibility of imaging with an auxiliary {\it passive} array, 
with sensors located at $( \vec\bx_q )_{q=1}^{\Nq}$, and with the scattering
medium in between it and the surface source-receiver array
(see Figure \ref{fig:intro0}, right).
The signals recorded by the auxiliary array form the data matrix
\begin{equation}
\label{data}
\big\{ p(t,\vec\bx_q;\vec\bx_s), \, t \in \RR, \, q=1,\ldots,\Nq,\, s=1,\ldots,\Ns\big\} .
\end{equation}
How can we use the auxiliary passive array to get an image, and will it help in mitigating the effects of the strong scattering
between it and the active source-receiver array? If we think of the strong scattering as producing
signals that appear to come from spatially dispersed, statistically independent noisy sources, then we are in a daylight imaging
setup with ambient noise sources, which was analyzed in \cite{cross}. 
Daylight imaging here means illumination from behind the receiver array.
By analogy with  having illumination from $\Ns$
uncorrelated point sources at $(\vec\bx_s)_{ s=1}^{\Ns} $,
we expect that, even in the case of active impulsive sources,  the matrix of cross correlations
at the auxiliary array
\begin{equation}
\label{def:cross parax0}
{\mathcal C}  \big(\tau, \vec\bx_{q},\vec\bx_{q'} \big) = 
\sum_{s=1}^{\Ns}  \int_\RR p(t,\vec\bx_q;\vec\bx_s)p(t+\tau,\vec\bx_{q'};\vec\bx_s)dt  
\end{equation}
will behave roughly as if it is the full active array response matrix of the auxiliary array. 
This means that it can be used for
imaging using Kirchhoff migration:
\begin{equation}
\label{def:imag}
{\mathcal I}_{\rm C}  (\vec\by^S) = \frac{1}{\Nq^2}  \sum_{q,q'=1}^{\Nq} {\mathcal C}  \Big(
  \frac{ | \vec\bx_q -  \vec\by^S | + | \vec\by^S -\vec\bx_{q'} |}{c_0} , \vec\bx_q,\vec\bx_{q'} \Big)  .
\end{equation}
In our previous  paper \cite{virtual12} we showed analytically the striking result that 
the resolution of the images obtained this way is given by the Rayleigh resolution
formulas as if the medium was homogeneous \cite{resolution10}. 
This result is not obvious, simply because scattering by inhomogeneities
of waves emitted by the impulsive sources at the surface produces signals on the auxiliary passive array
that are very different from those coming from spatially uncorrelated noise sources.  
This is especially true when addressing randomly layered media, which do not provide any lateral diversity enhancement 
since scattering does not change the transverse wave vectors.  Nevertheless, it has been
anticipated and observed in exploration geophysics contexts \cite{bakulin,schuster,wap10} 
that imaging with the cross correlations of the auxiliary
array is very effective and produces images that are nearly as good as in a homogeneous medium.
We showed this for randomly layered media in the asymptotic regime studied in detail in \cite{book}
and for isotropic random media in the paraxial regime \cite{tappert,uscinski77,ishimaru,hoopsolna,parax}.

The reason why this kind of imaging works so well is because, by wave field reciprocity, the cross correlations ${\mathcal C}(\tau,\vec\bx_q,\vec\bx_{q'})$
can be given a time-reversal interpretation, a well known \cite{derode95} observation, 
and then we can apply the analysis of time reversal refocusing in
random media \cite{blomgren,book}. The main results in time reversal in random media are (i) the enhanced refocusing, which
means that the primary peak in the cross correlation will be observable only when the sensors in the auxiliary array are close enough
with time lag close to zero, and then Kirchhoff migration to image the reflector will not be affected, and (ii) the statistical stability
of the cross correlation function relative to the random medium inhomogeneities, provided that the source illumination from the surface is
broadband (see \cite[Chapters 12 and 15]{book} and \cite{lord,fouque06}).

\begin{figure}
\begin{center}
\begin{tabular}{c}
\includegraphics[width=5.5cm]{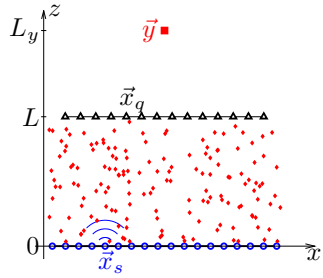} 
\end{tabular}
\end{center}
\caption{Sensor array imaging of a reflector  through a scattering medium in the region $z \in (0,L)$. $\vec\bx_s$ is a source,
$\vec\bx_q$ is a receiver, and $\vec\by$ is the reflector.}
\label{fig:intro1}
\end{figure}

In this paper we explore analytically the effect of finite array sizes in correlation-based imaging.
We first formulate the direct scattering problem more precisely.
The space coordinates are denoted by $\vec\bx=(\bx,z)\in \RR^2 \times \RR$.
The waves are emitted by a point source located at $\vec\bx_s$ which belongs to an
 array of sources $(\vec\bx_s)_{s=1,\ldots,\Ns}$ located in the plane 
$z=0$. The waves are recorded by an array of receivers $(\vec\bx_q)_{q=1,\ldots,\Nq}$  located 
in the plane $z=L$ (see Figure \ref{fig:intro1}). The recorded signals form the data matrix (\ref{data}).
The scalar wave field $(t,\vec\bx) \mapsto p(t,\vec\bx;\vec\bx_s)$ satisfies the wave equation
\begin{equation}
\label{eq:wave0}
\frac{1}{c (\vec\bx)^2} 
\frac{\partial^2 p}{\partial t^2} - \Delta p = - \nabla \cdot   \vec{\itbf F}  (t,\vec\bx;\vec\bx_s),
\end{equation}
where $c(\vec\bx)$ is the speed of propagation in the medium
and the forcing term $(t,\vec\bx) \mapsto  \vec{\itbf F}(t,\vec\bx;\vec\bx_s)$ models the source.
It is point-like,
located at $\vec\bx_s =(\bx_s,0^+)$, i.e. just above the surface $z=0$, and 
it emits a pulse:
\begin{equation}
\label{def:source}
\vec{\itbf F}  (t,\vec\bx;\vec\bx_s) =
\vec{\itbf f}  (t ) \delta(z) \delta( \bx -\bx_s )  ,
\end{equation}
We consider in this paper a randomly scattering medium that
occupies the section $z \in (0,L)$ and is sandwiched between two homogeneous half-spaces:
\begin{equation}
\frac{1}{c (\vec\bx)^2} = \frac{1}{c_0^2}\big( 1 + \mu(\vec\bx)\big),
\quad \quad \vec\bx \in \RR^2 \times(0,L),
\end{equation}
where $\mu(\vec\bx)$ is a zero-mean stationary random process modeling the random
heterogeneities present in the medium.
The homogeneous half-space $z >L$  is matched to the random section $z \in (0,L)$.

We consider scattering by a reflector above the random medium
 placed at $\vec\by =(\by,L_y)$, $L_y> L$.
The reflector is modeled by a local change of the speed of propagation of the form
\begin{equation}
\frac{1}{c (\vec\bx)^2} = \frac{1}{c_0^2}\Big( 1 + \frac{\sigma_{\rm ref}}{|\Omega_{\rm ref}|} {\bf 1}_{\Omega_{\rm ref}}(\vec\bx-\vec\by)\Big),
\quad \quad \vec\bx \in \RR^2 \times(L,\infty),
\end{equation}
where $\Omega_{\rm ref}$ is a small domain  
and $\sigma_{\rm ref}$ is the reflectivity of the reflector.\\
Finally we assume  Dirichlet boundary conditions at the surface $z=0$:
$$
p (t,(\bx,z=0))=0, \quad \quad \bx \in \RR^2.
$$

The recorded signals form the data matrix (\ref{data}).
The goal is to extract the location of the reflector from the data.
We will study the imaging function (\ref{def:imag}), considered in \cite{virtual12}, that migrates the 
cross correlation of the recorded signals (\ref{def:cross parax0}).
In \cite{virtual12} we analyzed the case in which the source array has full aperture and 
extends over the whole surface $z=0$. 
In this case we have shown both in the weakly scattering paraxial regime and in strongly scattering layered media
that the correlation-based imaging function (\ref{def:imag}) produces images as
if the medium between the sources and the receiver array was homogeneous
and the receiver array was an active one made up of both  sources and receivers.
This imaging method turns out to be very efficient as it completely cancels
the effect of random scattering.

In this paper we address the finite source aperture case, in which the sources do
not extend over the whole surface $z=0$. In this case it turns out that random scattering 
affects the resolution of the image, which is not the same
with and without random scattering.
However the effect of random scattering depends on its angular 
properties and it may enhance or reduce the angular diversity of the illumination,
which in turn enhances or reduces the resolution of the imaging function (\ref{def:imag}).
We will analyze in detail the weakly scattering paraxial regime, for which random
scattering is good for imaging,  and the strongly scattering layered regime,
for which random scattering is bad for imaging.

\section{Summary of the results}
 
We can give a simple explanation for why  the imaging function (\ref{def:imag})
gives a good image provided that some ideal conditions are fulfilled.
If the sources are point-like, generate Dirac-like pulses, 
and  densely surround the region of interest $\Omega$ (inside of which are the reflector
and the receiver array), then we have
$$
\hat{\mathcal C}(\omega,\vec\bx_q,\vec\bx_{q'}) =
 \int_{\partial \Omega}  \hat{G}(\omega,\vec\bx_q;\vec\bx_s) \overline{\hat{G}(\omega,\vec\bx_{q'};\vec\bx_s)}
  d\sigma(\vec\bx_s)  ,
$$
where $\hat{G}(\omega,\vec\bx_q;\vec\bx_s)$ is the full time-harmonic Green's function 
for the wave equation (\ref{eq:wave0}) including the reflector.
In this paper the Fourier transform of a function $f(t)$ is defined by
$$
\hat{f}(\omega)= \int f(t) e^{i\omega t} dt .
$$
By the Helmholtz-Kirchhoff identity \cite{cross,wap10}, we find that, provided $\Omega$ is a ball
with large radius:
$$
\hat{\mathcal C}(\omega,\vec\bx_q,\vec\bx_{q'}) 
= \frac{\omega}{c_0} {\rm Im} \big( \hat{G}(\omega,\vec\bx_q;\vec\bx_{q'}) \big) .
$$
This shows that the cross correlation of the signals at two receivers $\vec\bx_q$ and 
$\vec\bx_{q'}$ looks like the signal recorded at $\vec\bx_q$ when $\vec\bx_{q'}$ is a source.
Therefore, Kirchhoff Migration of the cross correlation matrix should give a good image.
The use of the Helmholtz-Kirchhoff identity gives the desired result, but it obscures the important role 
of scattering when the source array has limited aperture. We will show in this paper 
that it is not required to have full aperture source array 
to get a good image with the imaging function (\ref{def:imag}),
but this result requires a deeper mathematical analysis than the often-used  Helmholtz-Kirchhoff identity.

Let us now consider the situation shown in Figure \ref{fig:intro1},
when the source array has full aperture and covers the surface $z=0$, 
the angular illumination of the reflector is ultra-wide and the illumination cone covers the 
receiver array. This situation was analyzed in detail in \cite{virtual12}. 
In this case  the  correlation-based imaging function (\ref{def:imag}) 
completely cancels
the effect of random scattering and the results are equivalent 
whatever the form of the scattering medium. The cross-range resolution of the imaging function
is given by the classical Rayleigh resolution formula $\lambda_0 (L_y-L)/a$, where $a$
is the receiver array aperture. The range resolution is limited by the noise source bandwidth $B$
and given by $c_0/B$.

The results are quite different when the source array has limited aperture $b$.
In this case scattering turns out to play a critical role, as it may enhance or reduce the angular diversity
of illumination of the reflector. This was already noticed in time reversal  
\cite{derode95,blomgren,book}:
When waves emitted by a point source and recorded by an array are time-reversed and re-emitted into the medium, the time-reversed waves refocus at the original source location, and refocusing is 
enhanced in a scattering medium compared to a homogeneous one. This is because of the multipathing induced by scattering which enhances the refocusing cone.
However this is the first time in which this result is clearly seen in an imaging context,
in which the backpropagation step is carried out numerically in a fictitious homogeneous medium,
and not in the physical medium. This requires the 
backpropagation of the cross-correlations of the recorded
signals, and not the signals themselves.

\begin{figure}
\begin{center}
\begin{tabular}{c}
\includegraphics[width=6.3cm]{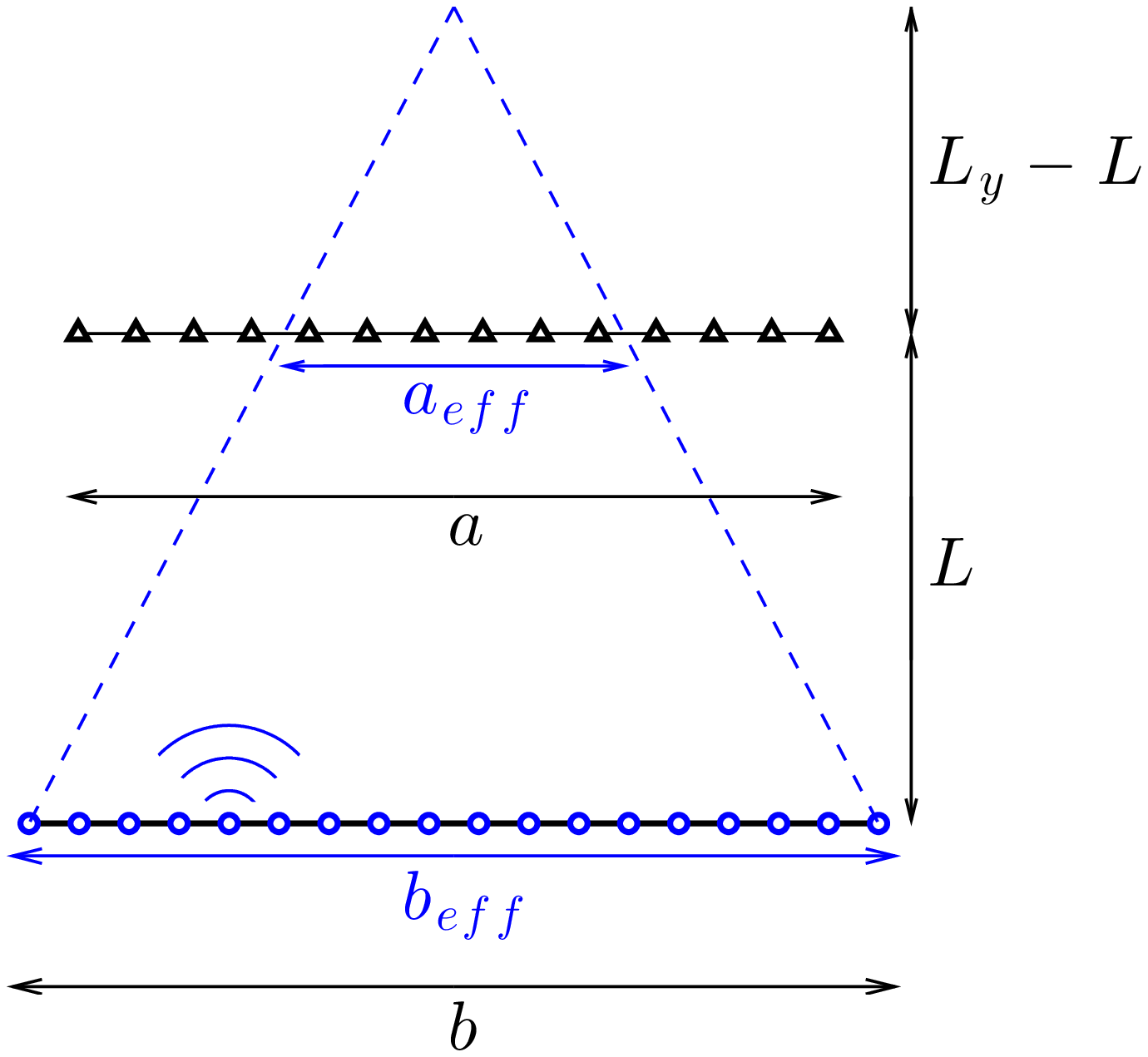} 
\hspace*{-0.5cm}
\includegraphics[width=6.3cm]{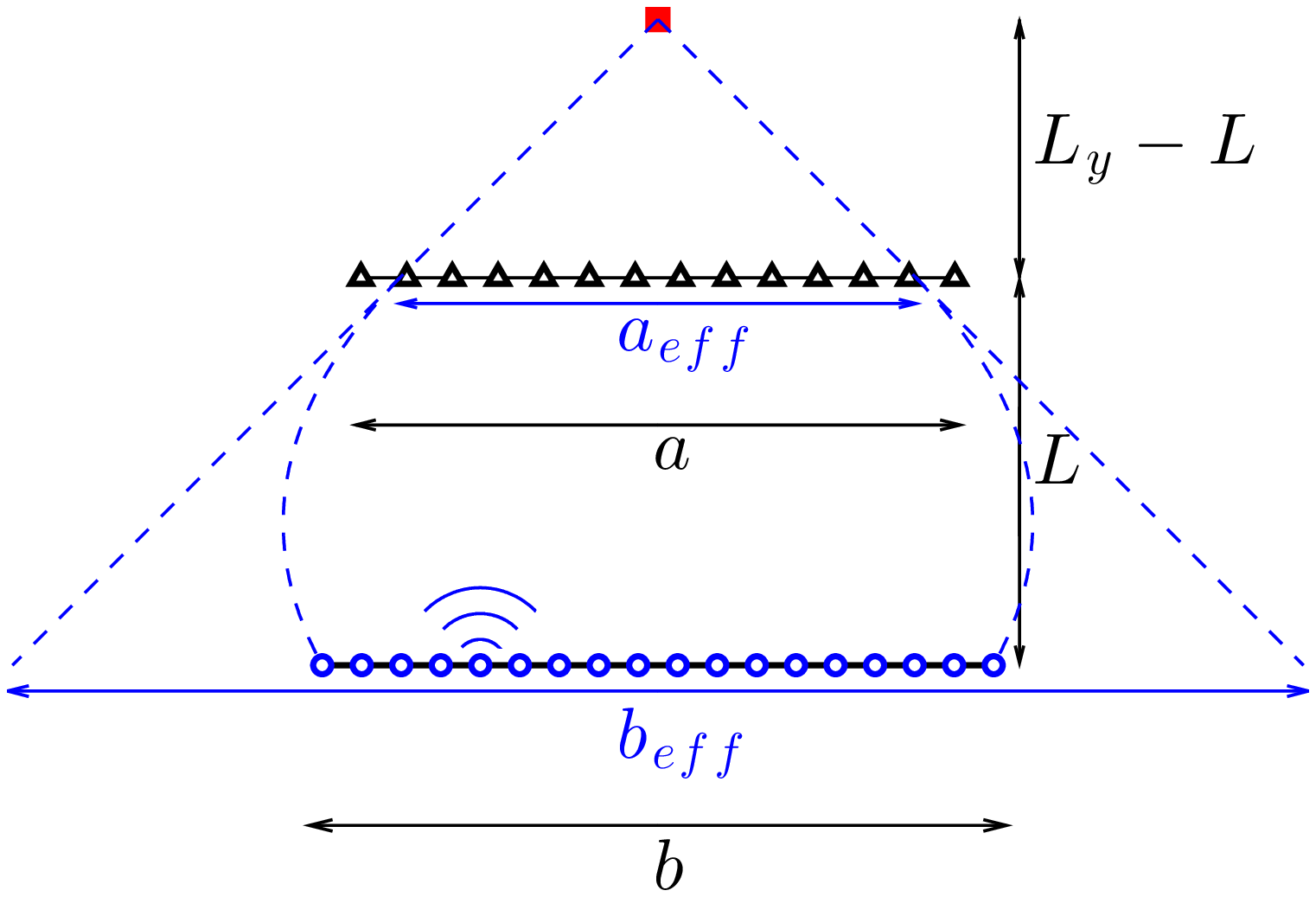} 
\end{tabular}
\end{center}
\caption{If the medium is homogeneous (left picture), the illumination aperture is the physical 
source aperture $b_{\rm eff}=b$. In  the random paraxial regime (right picture),
 the angular diversity of the 
waves that illuminate the reflector is increased by scattering and  the effective illumination 
aperture is enhanced $b_{\rm eff} > b$.}
\label{fig1}
\end{figure}

We first address the case of a medium with isotropic three-dimensional weak fluctuations $\mu(\vec\bx)$ of
the index of refraction. When the conditions for the paraxial approximation are fulfilled,
backscattering can be neglected and wave propagation is governed by a Schr\"odinger-type
equation with a random potential that has the form of a Gaussian field whose covariance 
function is given by
\begin{equation}
\label{def:covgaus}
\EE\big[ B(\bx,z)B(\bx',z')\big] =  \gamma_0(\bx-\bx') \big( z \wedge z' \big),
\end{equation}
with
\begin{equation}
\label{def:gamma0}
\gamma_0(\bx)=  \int_{-\infty}^\infty  \EE[ \mu({\bf 0},0)\mu(\bx,z) ] dz .
\end{equation}
We will show in Section \ref{app:paraxial}, by using multiscale analysis,
that the cone of incoherent waves that illuminates the reflector
is enhanced compared to the cone of coherent waves that illuminates the reflector through a homogeneous medium (see Figure \ref{fig1}), and this angular cone corresponds to an effective source
array aperture ${b}_{\rm eff}$ given by
\begin{equation}
\label{def:beff}
{b}_{\rm eff}^2 = {b}^2 +  \frac{\bar{\gamma}_2L^3 }{3}  ,
\end{equation}
where we have assumed that the correlation function $\gamma_0$ can be expanded as 
$\gamma_0(\bx) = \gamma_0({\bf 0}) - \bar{\gamma}_2 |\bx|^2+ o(|\bx|^2)$ for $|\bx|\ll 1$.
This in turn corresponds to an effective receiver array aperture ${a}_{\rm eff}$
(defined as the intersection of the illumination cone with the receiver array) given by:
\begin{equation}
\label{def:aeff}
 a_{\rm eff}= {b}_{\rm eff} \frac{L_y-L}{L_y}  .
\end{equation}
As a result, the  cross-range resolution of the imaging function
is given by the effective Rayleigh resolution formula $\lambda_0 (L_y-L)/a_{\rm eff}$,
which exhibits a resolution enhancement since $a_{\rm eff} $ is larger in a random medium
than in a homogeneous one.
The range resolution is still given by $c_0/B$.
The detailed analysis is in Subsection \ref{subsec:reso1}.

\begin{figure}
\begin{center}
\begin{tabular}{c}
\includegraphics[width=6.3cm]{fig_radar_1.eps} 
\hspace*{-0.5cm}
\includegraphics[width=6.3cm]{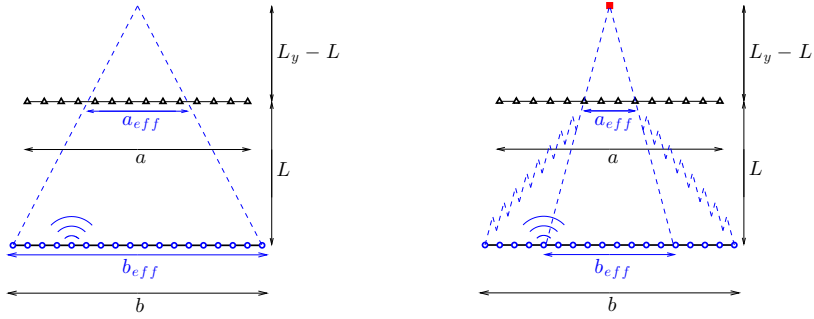} 
\end{tabular}
\end{center}
\caption{If the medium is homogeneous (left picture), the illumination aperture is the physical 
source aperture $b_{\rm eff}=b$. If the medium is randomly layered (right picture),
 the angular diversity of the 
waves that illuminate the reflector is reduced by scattering (only the waves with wavevectors
close to the vectical direction can reach the reflector after multiple scattering events
that conserve the wavevectors in the layered case)
and  the effective illumination 
aperture is reduced $b_{\rm eff} < b$.}
\label{fig2}
\end{figure}

We next analyze the case of a medium with one-dimensional (layered) fluctuations $\mu(z)$ of
the index of refraction.  In this case it is known \cite{book} that the scattering regime
is characterized by strong backscattering and wave localization, with the localization length given by:
\begin{equation}
\label{def:lloc}
L_{\rm loc} = \frac{4c_0^2}{\gamma \omega_0^2} ,
\end{equation}
where $\omega_0$ is the noise source central frequency and 
\begin{equation}
\label{def:gamma}
 \gamma = \int_{-\infty}^{\infty} \EE[ \nu(z') \nu(z'+z)] dz 
\end{equation}
is the integrated covariance of the fluctuations of the index of refraction.
We will show in Section \ref{sec:layered} that the cone of incoherent waves that illuminates the reflector
is reduced compared to the cone of coherent waves that illuminates the reflector through a homogeneous medium (see Figure \ref{fig2}), because scattering does not change the transverse wavevector,
and this angular cone corresponds to an effective source
array aperture ${b}_{\rm eff}$ given by
\begin{equation}
\label{def:beff:layered}
{b}_{\rm eff}^2 = \frac{4 L_y^2 L_{\rm loc}}{L} ,
\end{equation}
where we have assumed that $b^2 \gg L L_{\rm loc}$.
This in turn corresponds to an effective receiver array aperture ${a}_{\rm eff}$ given by:
\begin{equation}
\label{def:aeff:layered}
 a_{\rm eff}= {b}_{\rm eff} \frac{L_y-L}{L_y}  .
\end{equation}
As a result the  cross-range resolution of the imaging function
is given by the effective Rayleigh resolution formula $\lambda_0 (L_y-L)/a_{\rm eff}$,
which exhibits a resolution reduction since $a_{\rm eff}$ is smaller in a randomly layered medium
than in a homogeneous medium.
Furthermore, as wave scattering is strongly frequency-dependent,
the effective bandwidth is reduced as well 
\begin{equation}
B_{\rm eff}= \frac{B}{\sqrt{ 1+ \frac{B^2 L}{4 \omega_0^2 L_{\rm loc}}}},
\end{equation}
and the range resolution is  given by $c_0/B_{\rm eff}$.

The comparative analysis of the random paraxial regime and the randomly layered
regime clearly exhibits the role of scattering in correlation-based virtual source imaging.
With a full aperture source array, scattering plays no role as the illumination of the 
reflector is ultra-wide whatever the scattering regime.
When the source array is limited, if scattering is isotropic, then it 
enhances the angular diversity of the illumination of the reflector and the resolution
of the correlation-based imaging function.
If it is anisotropic, then it 
reduces the angular diversity of the illumination of the reflector and the resolution
of the correlation-based imaging function.

We note, however, that a large physical source array (and/or braodband sources)
is necessary to ensure the statistical stability of the imaging function. This has been already 
addressed in detail in different contexts in \cite{borcea11b,garnier13}.

\section{Imaging through a Random Medium in the Paraxial Regime}
\label{app:paraxial}%

\subsection{The Paraxial Scaling Regime}
In this section we analyze a scaling regime in which scattering is isotropic and
weak, which allows us to use the random paraxial wave model to describe the wave propagation
in the scattering region. In this approximation, backscattering is negligible but
there is significant lateral scattering as the wave advances, and over long
propagation distances. Even though they are weak, these effects accumulate
and can be a limiting factor in imaging and communications if not mitigated
in some way. Wave propagation in random media in the paraxial regime
has been used extensively in underwater sound propagation as well as in
the microwave and  optical contexts in the atmosphere 
\cite{uscinski77,tappert}.
%\cite{jakeman77,uscinski77,tappert,hudson80}.
We formulate the regime of paraxial wave propagation in
random media with a scaling of parameters that allows detailed
and effective mathematical analysis \cite{parax}. It is described as follows.

1) We assume that the correlation length $l_c$ of the medium
is much smaller than the typical
propagation distance $L$. 
We denote by $\eps^2$ the ratio between 
the correlation length and the typical propagation distance:
$$
\frac{l_c}{L} \sim \eps^2 .
$$

2) We assume that 
the transverse width of the source $R_0$ 
and the correlation length of the medium $l_c$ are of the same order.
This means that  the ratio $R_0/L$ is of order $\eps^2$.
This scaling is motivated by the fact that,
in this regime, there is a non-trivial interaction between the 
fluctuations of the medium and the wave.

3) We assume that the typical wavelength $\lambda$ is much smaller 
than the propagation distance $L$,
more precisely, we assume that the ratio $\lambda/L$ is of order $\eps^4$.
This high-frequency scaling is motivated by the following considerations. 
The Rayleigh length for a beam with initial width $R_0$ and  central wavelength $\lambda$ is 
of the order of $R_0^2 /\lambda$ when there is no random fluctuation.
The Rayleigh length is the distance from beam 
waist where the beam area is doubled by diffraction \cite{born}. 
In order to get a Rayleigh length of the order 
of the propagation distance $L$, the ratio $\lambda / L$ 
must be of order $\eps^4$ since $R_0/L \sim \eps^2$:
$$
\frac{\lambda}{L} \sim \eps^4   .
$$

4) We assume that the typical amplitude of the random fluctuations of the medium
is small. More precisely, we assume that the relative amplitude of the fluctuations is of order $\eps^3$.
This scaling has been chosen so as to obtain an effective regime of order one when $\eps$ goes to zero. 
That is, if the magnitude of the fluctuations is smaller than $\eps^3$, then the wave would 
propagate as if the medium was homogeneous, while if the order of magnitude is 
larger, then the wave would not be able to penetrate the random medium.
The scaling that we consider here corresponds to the physically most interesting 
situation where random effects play a role.

\subsection{The Random Paraxial Wave Equation}
We consider the time-harmonic form of the scalar wave equation
\begin{equation}
(\partial_z^2+\Delta_\perp) \hat{p} + \frac{\omega^2}{c_0^2} \big(1 + \mu(\bx,z)\big) \hat{p}=0  .
\end{equation}
Here $\mu$ is a zero-mean, stationary,  three-dimensional random process
with mixing properties in the $z$-direction.
In the high-frequency regime described above,
\begin{equation}
\label{eq:scaleparax}
\omega \to \frac{\omega}{\eps^4} , 
\quad  \quad  
\mu(\bx,z)  \to \eps^3 \mu\big(  \frac{\bx}{\eps^2}, \frac{z}{\eps^2}\big)  ,
\end{equation}
the rescaled function $\hat{\phi}^\eps$ defined by
\begin{equation}
\label{eq:scaleparax1}
\hat{p}^\eps(\omega,\bx,z)  = \exp \Big(  i \frac{\omega}{\eps^4} \frac{z}{c_0} \Big) \hat{\phi}^\eps \Big(\frac{\omega}{\eps^4},\frac{\bx}{\eps^2},z \Big) 
\end{equation}
satisfies
\begin{equation}
\label{eq:scaleparax2}
\eps^4 \partial_{z}^2  \hat{\phi}^\eps
+ 
\left(  2 i \frac{\omega}{c_0} \partial_z \hat{\phi}^\eps+ \Delta_\perp \hat{\phi}^\eps + \frac{\omega^2}{\eps c_0^2}  \mu\big(  \bx  , \frac{z}{\eps^2}\big)  \hat{\phi}^\eps \right)=0  .
\end{equation}
The ansatz (\ref{eq:scaleparax1}) corresponds to an up-going plane wave 
with a slowly varying envelope.
In the regime $\eps \ll 1$, it has been shown in \cite{parax} that
the forward-scattering approximation in the negative $z$-direction and the white-noise approximation are valid,
which means that the second-order derivative in $z$ in (\ref{eq:scaleparax2}) can be neglected and 
the random potential $\frac{1}{\eps}  \mu\big(  \bx  , \frac{z}{\eps^2}\big)$ can be replaced by a white noise.
The mathematical statement is that  the function $\hat{\phi}^\eps(\omega,\bx,z)$ converges to the solution $\hat{\phi}(\omega,\bx,z)$ of the
It\^o-Schr\"odinger equation
$$
 2 i  \frac{\omega}{c_0}d \hat{\phi} (\omega,\bx,z) +\Delta_\perp \hat{\phi} (\omega,\bx,z) dz + \frac{\omega^2}{ c_0^2}   \hat{\phi}(\omega,\bx,z)\circ d B(\bx,z) =0  ,
$$
where $B(\bx,z)$ is a Brownian field, that is a Gaussian process with mean zero and covariance function
(\ref{def:covgaus}). Here the $\circ$ stands for the Stratonovich stochastic integral~\cite{parax}.
We introduce the fundamental solution $\hat{G} \big(\omega,(\bx,z),(\bx_0,z_0) \big)$, which is defined as the solution
of the equation in $(\bx,z)$ (for $z > z_0$):
\begin{equation}
 2i  \frac{\omega}{c_0} d \hat{G} + \Delta_\perp \hat{G}  dz+  \frac{\omega^2}{c_0^2} \hat{G} \circ dB(\bx,z)= 0, 
\end{equation}
starting from $\hat{G}\big(\omega,(\bx,z=z_0),(\bx_0,z_0 ) \big) = \delta(\bx-\bx_0)$.
In a homogeneous medium ($B \equiv 0$) the fundamental solution is (for $z> z_0$)
\begin{equation}
\label{eq:green0}
\hat{G}_0 \big( \omega, (\bx,z), (\bx_0,z_0) \big) = \frac{  \omega}{2 i \pi c_0 (z-z_0)}   \exp \Big(  i \frac{\omega 
      |\bx-\bx_0|^2}{2c_0 (z-z_0)} \Big)    .
\end{equation}
In a random medium, the 
first two moments of the random fundamental solution have the following expressions.

\begin{proposition}
\label{prop:parax2}%
The first order-moment of the random  fundamental solution exhibits frequency-dependent damping:
\begin{equation}
\EE \big[ \hat{G}\big(\omega,(\bx,z),(\bx_0,z_0) \big)   \big] 
=
\hat{G}_0\big(\omega,(\bx,z),(\bx_0,z_0) \big) 
 \exp \Big( -\frac{\gamma_0({\bf 0}) \omega^2 |z-z_0|}{8 c_0^2}  \Big)  ,
\end{equation}
where $\gamma_0$ is given by (\ref{def:gamma0}).

The second order-moment of the random  fundamental solution exhibits spatial decorrelation:
\begin{eqnarray}
\nonumber
&& \EE \big[ \hat{G}\big(\omega,(\bx,z),(\bx_0,z_0) \big) 
\overline{\hat{G}\big(\omega,(\bx',z),(\bx_0,z_0) \big)} \big] \\
&& =
\hat{G}_0\big(\omega,(\bx,z),(\bx_0,z_0) \big) 
\overline{\hat{G}_0\big(\omega,(\bx',z),(\bx_0,z_0) \big)}
 \exp \Big( -  \frac{ \gamma_2(\bx-\bx') \omega^2 |z-z_0|}{4 c_0^2}   \Big)   , \hspace*{0.3in}
\end{eqnarray}
where $\gamma_2(\bx)= \int_0^1 \gamma_0({\bf 0}) -\gamma_0(\bx s) ds$.
\end{proposition}

These are classical results  \cite[Chapter 20]{ishimaru} once the the random paraxial equation has been proved to be correct,
as is the case here.
The result on the first-order moment shows that any coherent wave imaging method cannot give good images if the propagation distance is larger than 
the scattering mean free path $l_{\rm sca} = 8 c_0^2 /(\gamma_0({\bf 0}) \omega^2)$, because the coherent wave components will then be exponentially damped.
This is the situation we have in mind, and this is the situation in which imaging by migration of cross correlations turns out to be efficient.
The results on the second-order moment will be used in the next subsection to analyse the cross correlation of the recorded signals in a quantitative way.

\subsection{The Cross Correlation of Recorded Field in the Presence of a Reflector}
We consider the situation described in the introduction. In the random paraxial scaling regime
described above, the scalar field $p^\eps(t,\vec\bx;\vec\bx_s)$ corresponding to 
the emission from an element $\vec\bx_s$ of the surface source array
is solution of
\begin{equation}
\label{eq:wave01}
\frac{1}{c^\eps(\vec\bx)^2} 
\frac{\partial^2 p^\eps}{\partial t^2} - \Delta p^\eps = - \nabla \cdot  \vec{\itbf F}^\eps (t,\vec\bx;\vec\bx_s),
\end{equation}
where\\
- the source term is $\vec{\itbf F}^\eps (t,\bx,z;\vec\bx_s) =
\vec{\itbf f}^\eps  (t ) \delta(z) \delta( \bx -\bx_s )$,
the pulse 
is of the form
$$
\vec{\itbf f}^\eps (t) = f(t/\eps^4) \vec{\itbf e}_z  ,
$$
where the support of the Fourier transform 
of $f$ is bounded away from zero and of rapid decay
at infinity, and $\vec{\itbf e}_z$ is the unit vector pointing into the $z$-direction,\\
-  the medium is random in the region $z\in (0,L)$:
$$
\frac{1}{c^\eps(\vec\bx)^2 } = \frac{1}{c_0^2}\Big( 1 + \eps^3 \mu\big( \frac{\bx}{\eps^2} , \frac{z}{\eps^2}\big) \Big),
\quad \quad \vec\bx \in \RR^2 \times(0,L).
$$

We consider the cross correlation of the signals recorded at the receiver array  $(\vec\bx_q)_{q=1}^{\Nq}$ defined by:
\begin{equation}
\label{def:cross parax1}
{\mathcal C}^\eps \big(\tau, \vec\bx_{q},\vec\bx_{q'} \big) = 
\int_\RR \sum_{s=1}^{\Ns} p^\eps(t,\vec\bx_q;\vec\bx_s)p^\eps(t+\tau,\vec\bx_{q'};\vec\bx_s)dt  .
\end{equation}

Using the Born approximation for the point reflector at $\vec\by$, we obtain the following proposition \cite{virtual12}.

\begin{proposition}
In the random paraxial wave regime $\eps \to 0$, when there is a point reflector at $\vec\by=(\by,L_y)$ and when 
the source array covers the surface $z=0$, 
then the cross correlation of the recorded signals at the 
receiver array satisfies
\begin{eqnarray}
\nonumber
&&{\mathcal C}^\eps \Big( \frac{2L_y-2L}{c_0} + \eps^4 s, \vec\bx_{q},\vec\bx_{q'} \Big) \stackrel{\eps \to 0}{\longrightarrow}
 - \frac{i \sigma_{\rm ref}}{4\pi^3 c_0 (L_y-L)^2 } \int 
\omega^3  |\hat{f}(\omega)|^2 \\
&& \quad \quad \quad \quad \quad \quad \times \exp \Big( - i \omega\big( s -\frac{1}{2c_0} \frac{|\by-\bx_q|^2+ |\by-\bx_{q'}|^2}{L_y-L}\big) \Big) d\omega  .
\label{eq:cross parax2}
\end{eqnarray}
\end{proposition}

This shows that \\
- the cross correlation has a peak at time lag 
$$
\tau =  \frac{2L_y-2L}{c_0}  + \frac{\eps^4}{2c_0} \frac{|\by-\bx_q|^2+ |\by-\bx_{q'}|^2}{L_y-L},
$$
which is the sum of travel times from $\vec\bx_q$ to $\vec\by$ and from $\vec\by$ to $\vec\bx_{q'}$ in the paraxial
approximation:
\begin{eqnarray*}
\frac{|\vec\bx_q-\vec\by|}{c_0}+\frac{|\vec\by-\vec\bx_{q'}|}{c_0} &=&
\frac{1}{c_0} \sqrt{(L_y-L)^2 + \eps^4 |\by-\bx_q|^2} 
+
\frac{1}{c_0} \sqrt{(L_y-L)^2 + \eps^4 |\by-\bx_{q'}|^2} \\
&=& 
\frac{2L_y-2L}{c_0}
+
\frac{\eps^4}{2 c_0} \frac{|\by-\bx_q|^2+ |\by-\bx_{q'}|^2}{L_y-L} +O(\eps^8).
\end{eqnarray*}
- the effect of the random medium has completely disappeared.\\
The conclusion is therefore that Kirchhoff migration with cross correlations of the
receiver array produces images as if the medium was homogeneous and the 
receiver array was active.

When the source array has a finite aperture, say ${b}$, then an important quantity is the effective source array diameter
${b}_{\rm eff}$ defined by (\ref{def:beff})
The effective aperture can be interpreted as the source array diameter as seen from the receiver array through the random medium. 
It is increased by wave scattering induced by the random medium.
As we will see in the next section, 
this increase in turn enhances the resolution of the imaging function.

More precisely the following proposition shows that only the receivers 
that are within the cone determined by the effective source array  aperture contribute to the cross correlation.
As a result, the cross correlation is the same as in the case of a full aperture source array
provided that the effective array diameter is larger than a certain threshold value.
In the homogeneous case, this requires that the source array diameter must be larger than the 
threshold value. In the random medium case, the source array does not need to be large, only the 
effective source array diameter needs to be larger than the threshold value, which can be achieved 
 thanks to the second term in (\ref{def:beff}) which is due to scattering.

\begin{proposition}
We consider the random paraxial wave regime $\eps \to 0$, when there is a point reflector at $\vec\by=(\by,L_y)$ and when 
the source array covers a domain of radius ${b}$ at the surface $z=0$.
If the effective source array diameter is large enough in the sense that 
the effective Fresnel number $\frac{{b}_{\rm eff}^2}{\lambda_0 L} > \frac{L_y}{L_y-L}$, where $\lambda_0= 2\pi c_0/\omega_0$
 is the carrier wavelength,
then the cross correlation of the recorded signals at the 
receiver array satisfies
\begin{eqnarray}
\nonumber
&&{\mathcal C}^\eps \Big( \frac{2L_y-2L}{c_0} + \eps^4 s, \vec\bx_{q},\vec\bx_{q'} \Big) \stackrel{\eps \to 0}{Ê\longrightarrow}
 - \frac{i \sigma_{\rm ref}}{4\pi^3 c_0 (L_y-L)^2 } \int 
\omega^3  |\hat{f}(\omega)|^2 \psi_{\rm eff}( \bx_q ,\by) \\
&& \quad \quad \quad \quad \quad \quad \times \exp \Big( - i \omega\big( s -\frac{1}{2c_0} \frac{|\by-\bx_q|^2+ |\by-\bx_{q'}|^2}{L_y-L}\big) \Big) d\omega  .
\label{eq:cross parax2c}
\end{eqnarray}
with
\begin{equation}
 \psi_{\rm eff}(\bx_q , \by)= 
 \frac{{b}^2}{{b}_{\rm eff}^2}
  \exp \Big( - \frac{|\bx_q-\by L/L_y|^2}{a_{\rm eff}^2}
 \Big) .
\end{equation}
\end{proposition}

The finite source aperture array limits the angular diversity of the illumination, and as a result only a portion of the receiver
array really contributes to the cross correlation as characterized by the truncation function $\psi_{\rm eff}(\bx_q,\by)$.
In a homogeneous medium (left picture, figure \ref{fig1}) the truncation function has a clear geometric interpretation:
only the receivers localized along rays going from the sources to the reflector can contribute.
In a random medium, 
the angular diversity of the illumination is enhanced by scattering and
the truncation function is characterized by the effective source aperture ${b}_{\rm eff}$ 
that depends on the source array aperture ${b}$
and on the angular diversity enhancement induced by scattering (see (\ref{def:beff})).\\
$\bullet$
If $a_{\rm eff}  > a$, then the truncation function plays no role and we obtain the same result
as in the full aperture case.  \\
$\bullet$
If $a_{\rm eff}  < a$, then the truncation function plays a role and we obtain a result that is different from 
the full aperture case. \\
$\bullet$
In both cases scattering is helpful as it increases the angular diversity and reduces the impact of the truncation function.

\begin{proof}
We first describe the different wave signals that can be recorded at the surface array 
or at the receiver array.\\
1) The field recorded at the receiver passive array at $\vec\bx_q=(\bx_q,L)$ around time $L/c_0$ is the field transmitted through
the scattering layer in $z\in(0,L)$:
$$
p^\eps \big( \frac{L}{c_0} + \eps^4 s, \vec\bx_q;\vec\bx_s\big)  \stackrel{\eps \to 0}{\longrightarrow}
\frac{1}{2\pi} \int_{-\infty}^\infty \hat{f}(\omega) e^{-i \omega s} \hat{G} \big( \omega, (\bx_q,L), (\bx_s, 0) \big) d\omega   .
$$

2) Using the Born approximation for the reflector,  
the field recorded at the receiver passive array at $\vec\bx_{q'}=(\bx_{q'},L)$ around time $(2L_y-L)/c_0$ 
is of the form
\begin{eqnarray*}
p^\eps \big( \frac{2L_y-L}{c_0} + \eps^4 s, \vec\bx_{q'};\vec\bx_s \big)  \stackrel{\eps \to 0}{Ê\longrightarrow}
\frac{i\sigma_{\rm ref} c_0}{\pi} \int_{-\infty}^\infty \int_{\RR^2}
 \omega \hat{f}(\omega) e^{-i \omega s} \hat{G}_0 \big( \omega, (\by, L_y), (\bx_{q'},L) \big)\\
\times
 \hat{G}_0 \big( \omega,(\by, L_y) , (\bx,L) \big) \hat{G} \big( \omega,(\bx,L) , (\bx_s,0) \big) d\bx d\omega  .
\end{eqnarray*}
In the Born approximation there is no other wave component recorded at $\vec\bx_{q'}$ around a time $t_0 \not\in \{ L/c_0, (2L_y-L)/c_0\}$.\\
As a consequence the cross correlation of the signals recorded at the receiver array defined by
(\ref{def:cross parax1})
is concentrated around time lag $2(L_y-L)/c_0$ and it is of the form
\begin{eqnarray*}
{\mathcal C}^\eps \big( \frac{2L_y-2L}{c_0} + \eps^4 s, \vec\bx_{q},\vec\bx_{q'} \big)  \stackrel{\eps \to 0}{\longrightarrow}
\frac{i \sigma_{\rm ref} c_0}{\pi} \int_{-\infty}^\infty \int_{\RR^2}  \int_{\RR^2}
 \omega  |\hat{f}(\omega)|^2 e^{-i \omega s} \psi_{\rm s}(\bx_s)\\
 \times 
 \hat{G}_0 \big( \omega, (\by, L_y), (\bx_{q'},L) \big)   \hat{G}_0 \big( \omega,(\by, L_y) , (\bx,L) \big) \\
 \times
 \hat{G} \big( \omega,(\bx,L) , (\bx_s,0) \big) 
 \overline{\hat{G} \big( \omega, (\bx_q,L), (\bx_s, 0) \big)}d\bx_s d\bx d\omega  ,
\end{eqnarray*}
when the source array is dense at the surface $z=0$ and is characterized by the density function $\psi_{\rm s}$.
Using Proposition \ref{prop:parax2} and the self-averaging property of the product of two fundamental solutions
(one of them being complex conjugated)  \cite{PRS03,PRS07,parax} we get
\begin{eqnarray*}
{\mathcal C}^\eps \big( \frac{2L_y-2L}{c_0} + \eps^4 s, \vec\bx_{q},\vec\bx_{q'} \big)
 \stackrel{\eps \to 0}{Ê\longrightarrow}
\frac{i\sigma_{\rm ref}c_0}{\pi}   \int_{-\infty}^\infty \int_{\RR^2}  \int_{\RR^2}
 \omega  |\hat{f}(\omega)|^2 e^{-i \omega s} \psi_{\rm s}(\bx_s) \\
\times \hat{G}_0 \big( \omega, (\by, L_y), (\bx_{q'},L) \big)  
 \hat{G}_0 \big( \omega,(\by, L_y) , (\bx,L) \big) \\
 \times
 \hat{G}_0 \big( \omega,(\bx,L) , (\bx_s,0) \big) 
 \overline{\hat{G}_0 \big( \omega, (\bx_q,L), (\bx_s, 0) \big)} \\
 \times
  \exp \Big( -  \frac{\omega^2 L \gamma_2(\bx-\bx_q)}{4 c_0^2}   \Big)  
 d\bx_s d\bx d\omega  .
\end{eqnarray*}

When the sources cover the surface $z=0$, i.e. when $\psi_{\rm s} \equiv 1$,
we see by integrating in $\bx_s$ and by using the explicit expression (\ref{eq:green0})
that we have a Dirac distribution $\delta(\bx-\bx_q)$.
The exponential damping term then disappears because $\gamma_2({\bf 0})=0$ and
we finally obtain  (\ref{eq:cross parax2}).

When the source array has finite aperture with width ${b}$ 
at the surface $z=0$ and can be modeled by the density function
$\psi_{\rm s}(\bx_s)=\exp( -|\bx_s|^2/{b}^2)$, then
we get by integrating in $\bx_s$ and by using the explicit expression (\ref{eq:green0})
that 
\begin{eqnarray}
\nonumber
{\mathcal C}^\eps \big( \frac{2L_y-2L}{c_0} + \eps^4 s, \vec\bx_{q},\vec\bx_{q'} \big)
 \stackrel{\eps \to 0}{\longrightarrow}
\frac{\sigma_{\rm ref} {b}^2}{8\pi^3 c_0^2L^2(L_y-L)} 
\int  \omega^4 |\hat{f}(\omega)|^2  e^{-i \omega s}  \\
\times
 \hat{G}_0 \big( \omega,(\by, L_y) , (\bx_{q'},L) \big)
  {\mathcal G}( \omega, (\by, L_y), (\bx_{q},L) \big)
 d\omega  ,
 \label{eq:proof:parax2}
\end{eqnarray}
with
\begin{eqnarray*}
 {\mathcal G}( \omega, (\by, L_y), (\bx_{q},L) \big)&=&
 \int_{\RR^2}
 \exp \Big(  i \frac{\omega  |\by-\bx|^2}{2c_0 (L_y-L)}\Big)
 \exp \Big(  i \frac{\omega}{2c_0 L}( |\bx|^2 -|\bx_q|^2)\Big)
\\
&& \times   \exp \Big( -  \frac{\omega^2 {b}^2 |\bx-\bx_q|^2}{4 c_0^2 L^2}   \Big)  
  \exp \Big( -  \frac{\omega^2 L \gamma_2(\bx-\bx_q)}{4 c_0^2}   \Big)  
 d\bx
 \end{eqnarray*}

$\bullet$
 When there is no scattering or when scattering is weak in the sense that $\gamma_0({\bf 0}) \omega^2 L/c_0^2<1$
 for $\omega$ in the source bandwidth,
then  we have $\gamma_2\simeq  0$ and 
 $$
 \exp \Big( -  \frac{\omega^2 {b}^2 |\bx-\bx_q|^2}{4 c_0^2 L^2}   \Big)  
  \exp \Big( -  \frac{\omega^2 L \gamma_2(\bx-\bx_q)}{4 c_0^2}   \Big)  
\simeq 
 \exp \Big( -  \frac{\omega^2 {b}_{\rm eff}^2 |\bx-\bx_q|^2}{4 c_0^2 L^2}   \Big)  
$$
with ${b}_{\rm eff}={b}$.

$\bullet$
When scattering is strong so that $\gamma_0({\bf 0}) \omega^2 L/c_0^2>1$
for $\omega$ in the source bandwidth, then we have
$$
  \exp \Big( -  \frac{\omega^2 L \gamma_2(\bx-\bx_q)}{4 c_0^2}   \Big)  
  \simeq \exp \Big( -\frac{ \bar{\gamma}_2 \omega^2 L}{12 c_0^2} |\bx-\bx_q|^2 \Big)
$$
where $\bar{\gamma}_2$ such that $\gamma_0(\bx) = \gamma_0({\bf 0}) - \bar{\gamma}_2 |\bx|^2$ for $|\bx|\ll 1$,
so that $\gamma_2(\bx) = \bar{\gamma}_2 |\bx|^2/3$  for $|\bx|\ll 1$,
and
$$
 \exp \Big( -  \frac{\omega^2 {b}^2 |\bx-\bx_q|^2}{4 c_0^2 L^2}   \Big)  
  \exp \Big( -  \frac{\omega^2 L \gamma_2(\bx-\bx_q)}{4 c_0^2}   \Big)  
\simeq
 \exp \Big( -  \frac{\omega^2 {b}_{\rm eff}^2 |\bx-\bx_q|^2}{4 c_0^2 L^2}   \Big)  
$$
with ${b}_{\rm eff}^2={b}^2+\frac{\bar{\gamma}_2 L^3}{3}$.

By integrating in $\bx$ the expression of the function ${\mathcal G}$
we obtain 
\begin{eqnarray*}
 {\mathcal G}( \omega, (\by, L_y), (\bx_{q},L) \big)&=&
 \frac{\pi}{ \frac{\omega^2 {b}_{\rm eff}^2}{4 c_0^2 L^2} -\frac{i\omega L_y}{2c_0L(L_y-L)}}
\exp \Big\{ - \frac{ \frac{\omega^2}{4c_0^2}
  \big| \frac{\bx_q}{L}+\frac{\bx_q-\by}{L_y-L}\big|^2}{\frac{\omega^2 {b}_{\rm eff}^2}{4 c_0^2 L^2} -\frac{i\omega L_y}{2c_0L(L_y-L)}}  
  + \frac{i \omega |\bx_q-\by|^2}{2 c_0(L_y-L)} \Big\}
 \end{eqnarray*}
If $\frac{\omega_0 {b}_{\rm eff}}{c_0L} > \frac{L_y}{L_y-L}$, then
\begin{eqnarray*}
 {\mathcal G}( \omega, (\by, L_y), (\bx_{q},L) \big)&=&
 \frac{4 \pi c_0^2 L^2}{\omega^2 b_{\rm eff}^2}
\exp \Big\{ - \frac{  
  \big|  \bx_q - \frac{L}{L_y} \by\big|^2}{a_{\rm eff}^2}   
  + \frac{i \omega |\bx_q-\by|^2}{2 c_0(L_y-L)} \Big\} .
  \end{eqnarray*}
Substituting into (\ref{eq:proof:parax2}) gives the desired result.
\end{proof}

\subsection{Kirchhoff Migration of Cross Correlations} 
\label{subsec:reso1}%
The Kirchhoff migration function for the search point $\vec\by^S$ is
 \begin{equation}
 {\mathcal I}_{\rm C}^\eps (\vec\by^S) = \frac{1}{\Nq^2}  \sum_{q,q'=1}^{\Nq} {\mathcal C}^\eps \Big(
   \frac{ | \vec\bx_q -  \vec\by^S | + | \vec\by^S -\vec\bx_{q'} |}{c_0} , \vec\bx_q,\vec\bx_{q'} \Big)  ,
 \end{equation}
where $\Nq$ is the number of receivers at  the receiver array.
The following proposition describes the resolution properties of the imaging function
when the source array has full aperture. It was proved in \cite{virtual12}.

\begin{proposition}
If the receiver array at altitude $L$ is  a dense square array centered at ${\bf 0}$ and with sidelength $a$,
if the source array covers the surface $z=0$, 
if we assume  additionally  Hypothesis (\ref{hyp:H}):
\begin{equation}
\label{hyp:H}
\begin{array}{l}
\mbox{The sidelength $a$ is smaller than $L_y-L$.}\\
\mbox{The bandwidth $B$ of the source pulse is small compared}\\
\mbox{ to the central frequency~$\omega_0$.}
\end{array}
\end{equation}
Then, denoting the search point by
\begin{equation}
\label{searchpoint}
\vec\by^S= \vec\by + (\eps^2 \bxi, \eps^4\eta), 
\end{equation}
we have
 \begin{eqnarray}
 \nonumber
&& {\mathcal I}_{\rm C}^\eps (\vec\by^S) \stackrel{\eps \to 0}{\longrightarrow}
- \frac{i \sigma_{\rm ref}}{ 4\pi^3 c_0 (L_y-L)^2} 
 {\rm sinc}^2 \Big( \frac{\pi a \xi_1 }{ \lambda_0 (L_y-L)} \Big)
{\rm sinc}^2 \Big( \frac{\pi a \xi_2 }{ \lambda_0 (L_y-L)} \Big)
\\ 
&&\quad\quad  \quad \quad  \quad \quad\times
\int \omega^3 |\hat{f}(\omega)|^2 \exp \Big(  2 i \frac{\omega}{c_0} 
  \eta   \Big) d\omega   .
 \label{eq:funcIC:parax}
\end{eqnarray} 
\end{proposition}
Note that the result is not changed quantitatively if $a$ is of the same order as or larger than $L_y-L$
or if the bandwidth is of the same order as the central frequency, but then the transverse shape is not a sinc${}^2$ anymore.

This shows that the migration of the cross correlation gives the same result
as if we were migrating the array response matrix of the receiver array. 
Indeed  the imaging function (\ref{eq:funcIC:parax}) is exactly the imaging function that we would obtain 
if the medium was homogeneous, if the passive receiver array could be used as an active array, and if the 
 response matrix of the array was migrated to the search point $\vec\by^S$.
 In particular  the cross range resolution is $\lambda_0 (L-L_y) / a$
and the range resolution is $c_0/B$.

The following proposition is a new result
and it describes the resolution properties of the imaging function
when the source array has finite aperture.

\begin{proposition}
If the receiver array at altitude $L$ is  a dense square array centered at ${\bf 0}$ and with sidelength $a$,
if the source array has finite aperture ${b}$,
if we assume additionally Hypothesis (\ref{hyp:H}),
then, denoting the search point by  (\ref{searchpoint}) we have
 \begin{eqnarray}
 \nonumber
&& {\mathcal I}_{\rm C}^\eps (\vec\by^S) \stackrel{\eps \to 0}{\longrightarrow}
- \frac{i \sigma_{\rm ref}}{ 4\pi^3 c_0 (L_y-L)^2} 
 {\rm sinc} \Big( \frac{\pi a \xi_1 }{ \lambda_0 (L_y-L)} \Big)
{\rm sinc} \Big( \frac{\pi a \xi_2 }{ \lambda_0 (L_y-L)} \Big)
\\ 
 \nonumber
 &&\quad\quad  \quad \quad  \quad \quad\times
\frac{b^2}{a^2 b_{\rm eff}^2} \int_{[-a/2,a/2]^2} d\bx_q \exp \Big(-\frac{|\bx_q-\by|^2}{a_{\rm eff}^2} +
  i \frac{\omega_0}{c_0(L_y-L)}   \bxi\cdot (\bx_q-\by)  \Big)\\
&&\quad\quad  \quad \quad  \quad \quad\times
\int \omega^3 |\hat{f}(\omega)|^2 \exp \Big(  2 i \frac{\omega}{c_0} 
  \eta   \Big) d\omega   .
 \label{eq:funcIC:paraxfinite}
\end{eqnarray} 
\end{proposition}
This shows that:\\
- If the effective source aperture is large enough so that $|\bx_q - \by| \leq a_{\rm eff}$
for all $\bx_q \in [-a/2,a/2]^2$, then we get the same result (\ref{eq:funcIC:parax})
 as in the infinite source aperture case.\\
- If the effective source aperture is small,  then we get 
  \begin{eqnarray}
 \nonumber
&& {\mathcal I}_{\rm C}^\eps (\vec\by^S) \stackrel{\eps \to 0}{\longrightarrow}
- \frac{i \sigma_{\rm ref} b^2}{ 4\pi^2 c_0 L^2 a^2} 
 {\rm sinc} \Big( \frac{\pi a \xi_1 }{ \lambda_0 (L_y-L)} \Big)
{\rm sinc} \Big( \frac{\pi a \xi_2 }{ \lambda_0 (L_y-L)} \Big)
 \exp \Big(-\frac{\pi^2 a_{\rm eff}^2 |\bxi|^2}{\lambda_0^2(L_y-L)^2 }    \Big)\\
&&\quad\quad  \quad \quad  \quad \quad\times
\int \omega^3 |\hat{f}(\omega)|^2 \exp \Big(  2 i \frac{\omega}{c_0} 
  \eta   \Big) d\omega   ,
\end{eqnarray} 
which shows that the cross-range resolution is reduced (compared to the full source aperture case)
and the range resolution is not affected.

\section{Imaging through a Randomly Layered  Medium}
\label{sec:layered}%

In this section we analyze a scaling regime in which scattering is anisotropic and
strong. We  consider linear acoustic waves propagating in a three-dimensional
layered medium and generated by a point source.
Motivated by geophysical applications, we take a typical wavelength of the probing pulse
to be larger than the correlation length of the medium and smaller than the 
propagation distance. This is the regime appropriate in exploration geophysics studied for instance in \cite{asch91,book,garnier05,layered10b}.
In the analysis we abstract this regime of physical parameters by 
introducing dimensionless parameter $\eps$ that captures roughly the ordering of the scaling ratios:
\begin{equation}
\label{eq:scalingratio}
\frac{l_c}{L} \sim \eps^2, \quad \quad \frac{\lambda}{L} \sim \eps   . 
\end{equation}

We consider the situation described in the introduction. In the random paraxial scaling regime
described above, the scalar field $p^\eps(t,\vec\bx;\vec\bx_s)$ corresponding to 
an element $\vec\bx_s$ of the surface source array
is the solution of
\begin{equation}
\label{eq:wave02}
\frac{1}{c^\eps(\vec\bx)^2} 
\frac{\partial^2 p^\eps}{\partial t^2} - \Delta p^\eps = - \nabla \cdot  \vec{\itbf F}^\eps (t,\vec\bx;\vec\bx_s),
\end{equation}
where\\
- the source term is $\vec{\itbf F}^\eps (t,\vec\bx;\vec\bx_s) =
\vec{\itbf f}^\eps  (t ) \delta(z) \delta( \bx -\bx_s )$,
the pulse 
is of the form
\begin{equation}
\label{def:source2}
\vec{\itbf f}^\eps (t ) = \eps  \vec{\itbf f} \Big( \frac{t}{\eps} \Big)
,
\end{equation}
where we assume that the support of the Fourier transform 
of $\vec{\itbf f} =({\itbf f}_\bx,f_z)$ is bounded away from zero and of rapid decay
at infinity.
The particular scaling of $\vec{\itbf f}^\eps$ in (\ref{def:source})  means that the central wavelength is 
large compared to the microscopic scale  of variation of the random fluctuations of the medium and
small compared to 
the macroscopic scale of variation of the background medium, as in (\ref{eq:scalingratio}).
The normalizing amplitude factor $\eps$ multiplying the source term is not important 
as the wave equation is linear, but it will make the quantities of interest of order one as $\eps \to 0$,
which explains our choice.\\
-  the medium is randomly layered in the region $z\in (0,L)$:
\begin{equation}
\label{mod1a}
\frac{1}{c^\eps(\vec\bx)^2 } = \frac{1}{c_0^2}\Big( 1 + \nu(\frac{z}{\eps^2}) \Big),
\quad \quad \vec\bx=(\bx,z) \in \RR^2 \times(-L,0).
\end{equation}
In this model the parameters of the medium have 
 random and rapid fluctuations  with a typical scale of variation much smaller than 
the thickness of the layer.
The small dimensionless parameter $\eps^2$ is the ratio between these two scales.
The small-scale random fluctuations are  described by the random process $\nu(z)$.
The process $\nu$ is bounded in magnitude by a constant less than one, so that $c^\eps$ is  a positive quantity.
The random process $\nu(z)$ is stationary and zero mean.
It is supposed to have strong mixing 
properties so   that we can use averaging techniques for stochastic
differential equations as presented in \cite[Chapter 6]{book}.
The important quantity from the statistical point of view is 
the integrated covariance $\gamma$ of the fluctuations of the random medium defined by
(\ref{def:gamma}).
By the Wiener-Khintchine theorem it is nonnegative valued as it is the power spectral density evaluated at zero-frequency.
The integrated covariance $\gamma$ plays the role of 
scattering coefficient. It is of the order of the product of the variance of the medium fluctuations
times the correlation length of the medium fluctuations. As will become clear in the sequel, 
the statistics of the wave field depend on the random medium
via this integrated covariance or power spectral density.

%Motivated by geophysical applications we assume that the 
%density in the homogeneous half-space $z < 0$ is much smaller than the 
%typical density  in the medium for  $z \geq 0$.
%Since the velocity and pressure are continuous,
%the pressure in $z <  0$ goes to zero and hence, by continuity, also at $z=0$.
%These are the so-called pressure release boundary conditions:
%$$
%p (t,(\bx,z=0))=0, \quad \quad \bx \in \RR^2.
%$$

\subsection{Statistics of the Green's function}
When there is no reflector, the field is given by:
\begin{eqnarray}
\label{defp0}
 \hspace*{-0.3in}
  p^\eps (t,\vec\bx_q;\vec\bx_s) &=& -\frac{1}{(2\pi)^3 \eps } \int
{\mathcal G}^\eps_{\omega,\kappa} 
e^{-i \frac{\omega}{\eps} \left( t-\bk   \cdot   (\bx_q-\bx_s)-L / c_0(\kappa) \right)}
\hat{f}_z(\omega)
\omega^2 
d\omega  d\bk\,  .
\end{eqnarray}
Here

%$\bullet$
%The Fourier transforms are defined by 
%$$
%\hat{f}_z(\omega) = \int f_z(t) e^{i \omega t} dt.
%$$

$\bullet$
$c_0(\kappa)$ are the
 mode-dependent velocity
\begin{equation}
%\zeta_0(\kappa) =\frac{\zeta_0}{\sqrt{1-\kappa^2 c_0^2}} , \quad \quad
c_0(\kappa) =\frac{c_0}{\sqrt{1-\kappa^2 c_0^2}} . 
\end{equation}

$\bullet$
The random complex coefficient   ${\mathcal G}^\eps_{\omega,\kappa}$
is the Fourier-transformed Green's function (for Dirichlet boundary conditions at the surface $z=0$).
%pressure release boundary conditions).
The Fourier transform is taken both with respect to time and with respect to the transverse 
spatial variables.
The Green's function can be expressed in terms of the 
mode-dependent reflection and transmission coefficients 
 $R^\eps_{\omega,\kappa}$ and
$T^\eps_{\omega,\kappa}$ of the random section (for matched boundary conditions, 
that is, transparent boundary conditions) in the following way:
\begin{equation}
\label{generT}
{\mathcal G}^\eps_{\omega,\kappa} =
\frac{T^\eps_{\omega,\kappa}}{1- R^\eps_{\omega,\kappa}}  = \sum_{j=0}^\infty  
T^\eps_{\omega,\kappa} (R^\eps_{\omega,\kappa})^j 
 .
\end{equation}
When the medium is homogeneous, we have ${\mathcal G}^\eps_{\omega,\kappa}=1$.
When the medium is randomly layered, 
the statistics of $|{\mathcal G}^\eps_{\omega,\kappa}|^2$
was studied in \cite{layered10}. In particular it was shown that 
 $\EE[ |{\mathcal G}^\eps_{\omega,\kappa}|^2] \stackrel{\eps \to 0}{\longrightarrow} 1$,
 which is the result that is necessary and sufficient to study the cross correlation of the recorded signals
 for an infinite source aperture array.
 In the case of a finite source aperture array, the moment  of the Green's function at two nearby  frequencies  statistics is needed.
From \cite[Proposition 5.1]{layered10} we can obtain the second-order statistics of the Green's function.
\begin{proposition}
\label{prop:0}%
The autocovariance function of the Green's function at two nearby slownesses satisfies
\begin{equation}
\EE \big[   {\mathcal G}^\eps_{\omega,\kappa+\eps \lambda/2} \overline{{\mathcal G}^\eps_{\omega,\kappa-\eps \lambda/2}}
\big]
 \stackrel{\eps \to 0}{\longrightarrow}
\int {\mathcal U}(\omega,\kappa,\xi) \exp \big( -i \omega \kappa \lambda \xi \big) d\xi .
\end{equation}
The spectral density  ${\mathcal U}(\omega,\kappa,\xi)$ has a probabilistic representation.
For a fixed $(\omega,\kappa)$, it is the probability density function of a random variable 
\begin{equation}
\label{eq:represU}
{\mathcal U}(\omega,\kappa,\xi) =
\EE \Big[ \delta \Big(
 \xi-  2 c_0(\kappa)  \int_{0}^L N_{\omega,\kappa}(z)   dz \Big)
 \Big| N_{\omega,\kappa}(0) = 0 \Big]  
, 
\end{equation}
where $(N_{\omega,\kappa}(z))_{0 \leq z \leq L}$ is a jump Markov process with state space $\NN$ and inhomogeneous infinitesimal
generator
\begin{equation}
\label{def:l}
{\mathcal L}_z \phi(N) = \frac{\gamma c_0^2(\kappa ) \omega^2 }{4 c_0^4}
\big[ (N+1)^2 (\phi(N+1) -\phi(N)) + N^2( \phi(N-1)-  \phi(N))\big] .
\end{equation}
\end{proposition}

 The following proposition is new although it follows from the results obtained in \cite{book}.
 It characterizes the transition probabilities of the jump Markov process
 $(N_{\omega,\kappa}(z))_{0 \leq z \leq L}$ that is used in the probabilistic representation
 (\ref{eq:represU}) of the spectral density ${\cal U}$.

\begin{proposition}
\label{prop:0b}
The transition probabilities of the jump Markov process $N_{\omega,\kappa}$ are:
\begin{equation}
\PP \big( N_{\omega,\kappa}(z_0+z) = n \, \big| \, N_{\omega,\kappa}(z_0) = p \big)=
e^{- \frac{\widetilde{z}}{4}}
\int_0^\infty e^{-u^2 \widetilde{z}}Q_{n,p}(u) \frac{2 \pi u \sinh(\pi u)}{\cosh^2(\pi u)}
du ,
\end{equation}
where $\widetilde{z}= \frac{\gamma c_0^2(\kappa ) \omega^2 }{4 c_0^4} z$,
\begin{equation}
Q_{n,p}(u) = P_{n}(u)P_{p}(u) ,\quad \quad u \in [0,\infty), \quad n,p\in \NN,
\end{equation}
and
\begin{equation}
P_n(u)=\sum_{j=0}^n \binom{j}{n} (-1)^j K_{j}(u) ,\quad \quad
K_0(u)  =1,\quad \quad K_n(u) = \prod_{j=1}^{n} \frac{1}{j^2} \big[ u^2 +(j-\frac{1}{2})^2\big] ,
\end{equation}
or equivalently
\begin{eqnarray}
&&P_{n+1}(u) = \frac{2n^2+2n+\frac{3}{4}-u^2}{(n+1)^2} P_n(u)  -\frac{n^2}{(n+1)^2} P_{n-1}(u),\\
&&P_0(u)=1,\quad \quad P_1(u)=\frac{3}{4}-u^2 .
\end{eqnarray}
\end{proposition}
The polynomials $P_n(u)$ are orthonormal for the weight function 
$ \frac{2 \pi u \sinh(\pi u)}{\cosh^2(\pi u)}{\bf 1}_{[0,\infty)}(u)$.

\subsection{The Integral Representation of the Field}
In the presence of a reflector at $\vec\by=(\by,L_y)$, 
the scalar field at the position  $\vec\bx_q=(\bx_q,L)$  can be decomposed in a primary field and a secondary field:
$$
  p^\eps (t,\vec\bx_q;\vec\bx_s) =
  p_{\rm p}^\eps (t,\vec\bx_q;\vec\bx_s) 
  +  p_{\rm s}^\eps (t,\vec\bx_q;\vec\bx_s) .
$$
The primary field is the field obtained in the absence of the reflector given by:
\begin{eqnarray}
\label{defp0chap12}
 \hspace*{-0.3in}
  p_{\rm p}^\eps (t,\vec\bx_q;\vec\bx_s) &=& -\frac{1}{(2\pi)^3 \eps } \int
{\mathcal G}^\eps_{\omega,\kappa} 
e^{-i \frac{\omega}{\eps} \left( t-\bk   \cdot   (\bx_q-\bx_s)-L / c_0(\kappa) \right)}
\hat{f}_z(\omega)
\omega^2 
d\omega  d\bk\,  .
\end{eqnarray}

The secondary field is the additional contribution due to the reflector 
localized at $\vec\by=(\by,L_y)$ and given by
\begin{eqnarray}
\nonumber
 p_{\rm s}^\eps (t ,\vec\bx_q;\vec\bx_s) = \frac{i\sigma_{\rm ref}}{2(2\pi)^5 \eps^2}  
\int 
\Big(  \frac{1}{c_0(\kappa')} - c_0(\kappa) \bk   \cdot   \bk'
\Big)  
 \hat{f}_z(\omega)
\omega^5 {\mathcal G}^\eps_{\omega,\kappa'} 
\\
\times e^{ i \frac{\omega}{\eps} \phi(\bk,\bk')}
 d\omega  d \bk  d\bk' \, ,
\label{eq:defp1I}
\end{eqnarray}
where the  rapidly  varying phase is
\begin{eqnarray}
\label{eq:rapidphaseI}
\phi(\bk,\bk') &=& -t + \bk   \cdot  ( \bx_q-\by) + \frac{L_y-L}{c_0(\kappa)}
 +
\bk'   \cdot (  \by-\bx_s) + \frac{ L_y }{ c_0(\kappa')}\, .
\label{eq:rapidphaseII}
\end{eqnarray}
This expression has been obtained using the Born approximation for the reflector.

\subsection{The Cross Correlation of Recorded Field in the Presence of a Reflector}
The cross correlation ${\mathcal C}^\eps$ of the signals recorded at the receiver array is defined by
(\ref{def:cross parax1}). 
If the source array is dense with density $\psi_{\rm s}(\bx_s)$, the cross correlation
has the form 
\begin{equation}
{\mathcal C}^\eps(\tau ,\vec\bx_q,\vec\bx_{q'})  =
\eps^{-2} \int_{\RR^2} d\bx_s \psi_{\rm s}(\bx_s) \int dt p^\eps(t,\vec\bx_q;\vec\bx_s) p^\eps( t+\tau,\vec\bx_{q'};\vec\bx_s) ,
\end{equation}

In the presence of a reflector at $\vec\by$, it 
can be written as the sum 
\begin{equation}
{\mathcal C}^\eps(\tau ,\vec\bx_q,\vec\bx_{q'})  = 
{\mathcal C}_{\rm pp}^\eps(\tau ,\vec\bx_q,\vec\bx_{q'})  +
{\mathcal C}_{\rm ps}^\eps(\tau ,\vec\bx_q,\vec\bx_{q'})  +
{\mathcal C}_{\rm sp}^\eps(\tau ,\vec\bx_q,\vec\bx_{q'})  .
\end{equation}
Here ${\mathcal C}_{\rm pp}^\eps$ is the cross correlation of the primary field (\ref{defp0chap12})
 at $\vec\bx_q=(\bx_q,L)$ with the primary field at $\vec\bx_{q'}=(\bx_{q'},L)$,
\begin{equation}
{\mathcal C}_{\rm pp}^\eps(\tau ,\vec\bx_q,\vec\bx_{q'})  =
\eps^{-2} \int_{\RR^2} d\bx_s \psi_{\rm s}(\bx_s) \int dt p_{\rm p}^\eps(t,\vec\bx_q;\vec\bx_s) p_{\rm p}^\eps( t+\tau,\vec\bx_{q'};\vec\bx_s) ,
\end{equation}
${\mathcal C}_{\rm ps}^\eps$ is the cross correlation of the primary field  (\ref{defp0chap12}) 
at $\vec\bx_q$ with the secondary field (\ref{eq:defp1I})
 at $\vec\bx_{q'}$, 
\begin{equation}
{\mathcal C}_{\rm ps}^\eps(\tau ,\vec\bx_q,\vec\bx_{q'})  =
\eps^{-2} \int_{\RR^2} d\bx_s \psi_{\rm s}(\bx_s)\int dt p_{\rm p}^\eps(t,\vec\bx_q;\vec\bx_s) p_{\rm s}^\eps( t+\tau,\vec\bx_{q'};\vec\bx_s) ,
\end{equation}
and  ${\mathcal C}_{\rm sp}^\eps$ is
 the cross correlation of the secondary field (\ref{eq:defp1I}) at $\vec\bx_q$ with the primary field
 (\ref{defp0chap12}) at $\vec\bx_{q'}$,
\begin{equation}
{\mathcal C}_{\rm sp}^\eps(\tau ,\vec\bx_q,\vec\bx_{q'})  =
\eps^{-2} \int_{\RR^2} d\bx_s \psi_{\rm s}(\bx_s)\int dt p_{\rm s}^\eps(t,\vec\bx_q;\vec\bx_s) p_{\rm p}^\eps( t+\tau,\vec\bx_{q'};\vec\bx_s) .
\end{equation}

\begin{proposition}
\label{prop:crossps}%
For $\vec\by=(\by,L_y)$ and $\vec\bx_q=(\bx_q,L)$, 
the $ps$-secondary cross correlation centered at 
$\tau = [ |  \vec\by - \vec\bx_q | + | \vec\bx_{q'} - \vec\by | ]/c_0$ has the form:
\begin{eqnarray}
\nonumber
&& {\mathcal C}_{\rm ps}^\eps
\Big(  \frac{ |  \vec\by - \vec\bx_q | + | \vec\bx_{q'} - \vec\by |}{c_0}+\eps s  ,\vec\bx_q,\vec\bx_{q'} \Big) \\
\nonumber
 &&  \stackrel{\eps \to 0}{\longrightarrow}
\frac{i \sigma_{\rm ref}}{2 (2\pi)^3 c_0^3} \frac{|L_y-L|   
 \, (\vec\by - \vec\bx_q )\cdot (\vec\bx_{q'}-\vec\by ) }
{ |  \vec\by - \vec\bx_q |^3 \, |\vec\bx_{q'}- \vec\by  |^2} 
\\
\nonumber
&& \hspace*{0.5in} \times 
\int  \Psi_{\rm eff} \big( \omega,\vec\bx_q  , \vec\by \big)
  |\hat{f}_z(\omega)|^2 \omega^3 \exp\big( i \omega s\big) d\omega , \quad 
\end{eqnarray}
where
\begin{eqnarray}
\nonumber
 \Psi_{\rm eff} \big( \omega, \vec\bx_q  ,\vec\by  \big) =  \int_\RR \psi_{\rm s} \Big( \bx_q - (\by-\bx_q)
\big( \frac{\xi}{c_0 |\vec\by-\vec\bx_q|} +\frac{L}{L_y-L}\big) \Big) \\
\times
{\mathcal U} \Big(\omega,\frac{|\by-\bx_q|}{c_0 |\vec\by-\vec\bx_q|} ,\xi\Big) d\xi  .
\end{eqnarray}
The $ps$-secondary cross correlation is negligible elsewhere.
\end{proposition}

Let us first consider the case in which there is no scattering. Then the spectral density is
$$
{\mathcal U}  (\omega,\kappa,\xi )=\delta(\xi)
$$
and we have
\begin{eqnarray}
\nonumber
&& {\mathcal C}_{\rm ps}^\eps
\Big(  \frac{ |  \vec\by - \vec\bx_q | + | \vec\bx_{q'} - \vec\by |}{c_0}+\eps s  ,\vec\bx_q,\vec\bx_{q'} \Big) \\
\nonumber &&  \stackrel{\eps \to 0}{\longrightarrow}
\frac{i \sigma_{\rm ref}}{2 (2\pi)^3 c_0^3} 
\frac{|L_y-L|   
 \, (\vec\by - \vec\bx_q )\cdot (\vec\bx_{q'}-\vec\by ) }
{ |  \vec\by - \vec\bx_q |^3 \, |\vec\bx_{q'}- \vec\by  |^2} \psi_{\rm s} \big( {\itbf X}(\bx_q) \big)
%\\
%&& \hspace*{1.5in} \times 
\int |\hat{f}_z(\omega)|^2 \omega^3 \exp\big( i \omega s\big) d\omega  , \quad 
\end{eqnarray}
where
$$
{\itbf X}(\bx_q)  =  \bx_q- (\by-\bx_q) \frac{L}{L_y-L} 
$$
is the intersection of the ray going through $\vec\by$ and $\vec\bx_q$ with the surface $z=0$.
The geometric interpretation is clear: in absence of scattering, the receiver at $\vec\bx_q$ can contribute to the
cross correlation only if there is a ray going from the source to the target through it.
If $\psi(\bx)={\bf 1}_{[-b/2,b/2]^2}(\bx)$ and $\by={\bf 0}$, then we find that
$\psi({\itbf X}(\bx_q))={\bf 1}_{[-a_{\rm eff}/2,a_{\rm eff}/2]^2}(\bx_q)$
with $a_{\rm eff}= (L_y-L) b/L_y$ (see left picture, Figure \ref{fig2}).
If $a_{\rm eff} <a$, then this means that the receiver array aperture cannot be used 
to its maximal capacity because of limited source illumination.

Let us consider the case of a randomly layered medium.
In the strongly scattering regime the truncation function $\Psi_{\rm eff}$ takes a simple form.

\begin{lemma}
When there is strong scattering in the sense that $L \gg L_{\rm loc}$ and when $\psi$ is compactly supported 
in some bounded domain with diameter $b \ll \sqrt{L L_{\rm loc}}$, then 
$$
 \Psi_{\rm eff} \big( \omega, (\bx_q ,L) ,(\by,L_y) \big)\simeq 
C_\psi \Big( \frac{L_{\rm loc}}{L}\Big)^{3/2} \exp \Big(- \frac{L}{4L_{\rm loc} }  \frac{\omega^2}{\omega_0^2} 
 - \frac{L}{4L_{\rm loc} }  \frac{|\by-\bx_q|^2}{(L_y-L)^2}\Big) ,
$$
for some constant $C_\psi$ that depends on the source array aperture.
\end{lemma}

Note that the truncation function can also be written as
$$
 \Psi_{\rm eff} \big( \omega, (\bx_q ,L) ,(\by,L_y) \big)\simeq 
C_\psi \Big( \frac{L_{\rm loc}}{L}\Big)^{3/2} \exp \Big(- \frac{L}{4L_{\rm loc} } 
 \frac{\omega^2}{\omega_0^2} 
 -  \frac{|\by-\bx_q|^2}{a_{\rm eff}^2}\Big) ,
$$
in terms of the effective receiver array aperture $a_{\rm eff}$ defined by (\ref{def:beff:layered}-\ref{def:aeff:layered}).
This truncation function determines which receivers contribute to the imaging function
of the reflector located at $(\by,L_y)$.

\begin{proof}
We use the probabilistic representation of ${\mathcal U}$
given in Proposition \ref{prop:0}.
We note  first that, since $N_{\omega,\kappa}(z)$ takes nonnegative values, the function $\xi \to 
{\mathcal U} (\omega,\kappa,\xi) $ is supported in $[0,\infty)$.
From Propositions \ref{prop:0}-\ref{prop:0b}, when $\widetilde{z}\gg 1$ and $n,p = O(1)$, we have
\begin{equation}
\PP \big( N_{\omega,\kappa}(z_0+z) = n \, \big| \, N_{\omega,\kappa}(z_0) = p \big)=
\exp\Big(-  \frac{\widetilde{z}}{4}\Big)
\frac{\pi^{5/2} Z_{n,p}(0)}{2 \widetilde{z}^{3/2}} \big( 1 + o (1) \big) .
\end{equation}
As a consequence, for any $M$, there exists $C_M$ such that 
\begin{equation}
\PP \Big( \int_{0}^L N_{\omega,\kappa}(z) dz \leq M \Big)=
\exp \Big( - \frac{\widetilde{L}}{4}\Big)
\frac{C_M}{  \widetilde{L}^{3/2}} \big( 1 + o (1) \big) ,
\end{equation}
when $\widetilde{L}=L/L_{\omega,\kappa} \gg 1$, where
$$
L_{\omega,\kappa} = L_{\rm loc} \frac{c_0^2}{c_0(\kappa)^2} \frac{\omega_0^2}{\omega^2} .
$$
With $\kappa= |\by-\bx_q|/(|\vec\by-\vec\bx_q| c_0)$, we have
$$
\frac{1}{L_{\omega,\kappa}} = \frac{1}{L_{\rm loc}}  \frac{\omega^2}{\omega_0^2}\Big( 1 + \frac{|\by-\bx_q|^2}{(L_y-L)^2} \Big)  ,
$$
which gives the desired result.
\end{proof}

\subsection{Kirchhoff Migration of Cross Correlations} 
The Kirchhoff migration function for the search point $\vec\by^S$ is
 \begin{equation}
 {\mathcal I}_{\rm C}^\eps (\vec\by^S) = \frac{1}{\Nq^2} \sum_{q,q'=1}^{\Nq} {\mathcal C}^\eps \Big(
   \frac{ | \vec\bx_q -  \vec\by^S | + | \vec\by^S -\vec\bx_{q'} |}{c_0} , \vec\bx_q,\vec\bx_{q'} \Big)  ,
 \end{equation}
where $\Nq$ is the number of receivers at  the receiver array.

When the source array aperture is infinite, we find that the Kirchhoff migration function
gives an image that does not depend on scattering and is the same one as if the medium were homogeneous
\cite{virtual12}.

\begin{proposition}
If the receiver array at altitude $L$ is  a dense square array with sidelength $a$
(i.e. the distance between the sensors is smaller than or equal to half-a-central wavelength),
if the source array covers the surface $z=0$,
then, denoting the search point by
 \begin{equation}
 \label{searchpointlayer}
\vec\by^S= \vec\by +\eps (\bxi,\eta), 
 \end{equation}
we have
 \begin{eqnarray}
 \nonumber
&& {\mathcal I}_{\rm C}^\eps (\vec\by^S) \stackrel{\eps \to 0}{\longrightarrow}
 \frac{i \sigma_{\rm ref} }{2 (2\pi)^3 c_0^3 (L_y-L)^2} 
 {\rm sinc}^2 \Big( \frac{\pi a \xi_1 }{ \lambda_0 (L_y-L)} \Big)
{\rm sinc}^2 \Big( \frac{\pi a \xi_2 }{ \lambda_0 (L_y-L)} \Big)
\\ 
&&\quad\quad  \quad \quad  \quad \quad\times
\int \omega^3 |\hat{f}_z(\omega)|^2 \exp \Big( 2 i \frac{\omega}{c_0} 
 \big( \eta -  \frac{\bxi \cdot \by}{L_y-L} \big) \Big) d\omega   .
 \label{eq:funcIC}
\end{eqnarray} 
\end{proposition}

When the source array aperture is finite, Kirchhoff migration does not give the same
image in the presence and in the absence of random scattering. 
In the randomly layered regime addressed in this section, random scattering reduces the angular diversity of the 
illumination of the region of interest below the random medium. As a result
the image resolution is reduced as well as shown by the following proposition.

\begin{proposition}
If the receiver array at altitude $L$ is  a dense square array centered at $\bx_{\rm A}$ and with sidelength $a$
if the source array has finite aperture $b$ at the surface $z=0$,
and if we denore the search point by (\ref{searchpointlayer}), then 
we have
 \begin{eqnarray}
 \nonumber
&& \hspace*{-0.2in}
{\mathcal I}_{\rm C}^\eps (\vec\by^S) \stackrel{\eps \to 0}{\longrightarrow}
 \frac{i \sigma_{\rm ref} }{2 (2\pi)^3 c_0^3 (L_y-L)^2} 
 {\rm sinc} \Big( \frac{\pi a \xi_1 }{ \lambda_0 (L_y-L)} \Big)
{\rm sinc} \Big( \frac{\pi a \xi_2 }{ \lambda_0 (L_y-L)} \Big)
\\ 
 \nonumber
&&\quad\quad   \quad \quad\times
\int   d\omega  \omega^3 |\hat{f}_z(\omega)|^2 \exp \Big(  2 i \frac{\omega}{c_0} 
\big( \eta -  \frac{\bxi \cdot \by}{L_y-L} \big) \Big) \\
&&\quad\quad  \quad \quad\times
\frac{1}{a^2}  \iint_{[-a/2,a/2]^2}   \Psi_{\rm eff} \big( \omega, (\bx_q ,-L) ,(\by,-L_y) \big)
  \exp \Big(  i \frac{\omega}{c_0} 
  \frac{\bxi \cdot \bx_q}{L_y-L}  \Big)d\bx_q
  . \quad
 \end{eqnarray} 
Assume $\by={\bf 0}$ for simplicity.

If  scattering is weak ($L \ll L_{\rm loc}$)  then
 \begin{eqnarray}
 \nonumber
&&  \hspace*{-0.2in} {\mathcal I}_{\rm C}^\eps (\vec\by^S) \stackrel{\eps \to 0}{\longrightarrow}
 \frac{i \sigma_{\rm ref} a_{\rm eff}^2}{2 (2\pi)^3 c_0^3 (L_y-L)^2 a^2} 
 {\rm sinc} \Big( \frac{\pi a \xi_1 }{ \lambda_0 (L_y-L)} \Big)
{\rm sinc} \Big( \frac{\pi a \xi_2 }{ \lambda_0 (L_y-L)} \Big)
\\ 
&& \quad\times
 {\rm sinc} \Big( \frac{\pi a_{\rm eff} \xi_1 }{ \lambda_0 (L_y-L)} \Big)
{\rm sinc} \Big( \frac{\pi a_{\rm eff} \xi_2 }{ \lambda_0 (L_y-L)} \Big)
\int  \omega^3 |\hat{f}_z(\omega)|^2 \exp \Big( 2 i \frac{\omega}{c_0} \eta  \Big)
  d\omega   ,
\end{eqnarray} 
with $a_{\rm eff}= \min \big( (L_y-L)b / L_y , a)$.

If scattering is strong ($L \gg L_{\rm loc}$) then
 \begin{eqnarray}
 \nonumber
&& {\mathcal I}_{\rm C}^\eps (\vec\by^S) \stackrel{\eps \to 0}{\longrightarrow}
 \frac{i \sigma_{\rm ref} L_{\rm loc}}{ (2\pi)^2 c_0^3La^2} 
 {\rm sinc} \Big( \frac{\pi a \xi_1 }{ \lambda_0 (L_y-L)} \Big)
{\rm sinc} \Big( \frac{\pi a \xi_2 }{ \lambda_0 (L_y-L)} \Big)
\exp \Big( -  \frac{\pi^2 L_{\rm loc} |\bxi|^2}{L \lambda_0^2} \Big)
\\ 
&&\quad\quad  \quad \quad  \quad \quad\times
\int  \omega^3 |\hat{f}_z(\omega)|^2 \exp \Big( - \frac{L \omega^2}{4 L_{\rm loc} \omega_0^2} \Big)
\exp \Big( 2 i \frac{\omega}{c_0} \eta  \Big)
  d\omega   .
\end{eqnarray} 
\end{proposition}

If, in particular,
$$
\hat{f}(\omega)= \exp \Big( -\frac{(\omega-\omega_0)^2}{2B^2}\Big) ,
$$
then, in the weak  scattering regime ($L \ll L_{\rm loc}$) we have
 \begin{eqnarray}
 \nonumber
&& {\mathcal I}_{\rm C}^\eps (\vec\by^S) \stackrel{\eps \to 0}{\longrightarrow}
 \frac{i \omega_0^3 \sqrt{\pi} B  \sigma_{\rm ref} a_{\rm eff}^2}{2 (2\pi)^3 c_0^3 (L_y-L)^2 a^2} 
 {\rm sinc}^2 \Big( \frac{\pi a \xi_1 }{ \lambda_0 (L_y-L)} \Big)
{\rm sinc}^2 \Big( \frac{\pi a \xi_2 }{ \lambda_0 (L_y-L)} \Big)
\\ 
&&\quad\quad  \quad \quad  \quad \quad\times
\exp \Big( - \frac{B^2 \eta^2}{c_0^2}\Big),
\end{eqnarray} 
and in the strong  scattering regime ($L \gg L_{\rm loc}$) we have
 \begin{eqnarray}
 \nonumber
&& 
\hspace*{-0.2in}
{\mathcal I}_{\rm C}^\eps (\vec\by^S) \stackrel{\eps \to 0}{\longrightarrow}
 \frac{i \omega_0^3 \sqrt{\pi} B  \sigma_{\rm ref} L_{\rm loc}}{ (2\pi)^2 c_0^3La^2 \sqrt{1 + \frac{B^2 L}{4\omega_0^2 L_{\rm loc}}}} 
 \exp \Big( - \frac{L}{4 L_{\rm loc} \big(1 +  \frac{B^2 L}{4\omega_0^2 L_{\rm loc}}\big)}\Big) \\
 \nonumber
&&  \quad \quad  \quad \quad\times
 {\rm sinc} \Big( \frac{\pi a \xi_1 }{ \lambda_0 (L_y-L)} \Big)
{\rm sinc} \Big( \frac{\pi a \xi_2 }{ \lambda_0 (L_y-L)} \Big)
\exp \Big( -  \frac{\pi^2 L_{\rm loc} |\bxi|^2}{L \lambda_0^2} \Big)
\\ 
&&  \quad \quad  \quad \quad\times
\exp \Big( - \frac{B^2 \eta^2}{c_0^2 \big(1 +  \frac{B^2 L}{4\omega_0^2 L_{\rm loc}}\big)}\Big)
\exp \Big( - i \frac{B^2 L \eta}{2 \omega_0 L_{\rm loc} c_0 \big(1 +  \frac{B^2 L}{4\omega_0^2 L_{\rm loc}}\big)}
\Big) .
\end{eqnarray} 
We can see that both the cross-range resolution and the range resolution have 
been reduced by scattering. \\
The reduction in cross-range resolution comes the reduction in the effective illumination angular aperture
discussed above.\\
The reduction in range resolution comes from the fact that the high-frequency components
are more sensitive to scattering by thin random layers and therefore the effective spectrum
used in the imaging function is reduced compared to the original source spectrum.

\section*{Acknowledgments}
The work of G. Papanicolaou was
partially supported by AFOSR grant FA9550-11-1-0266.  The work of
J. Garnier was partially supported by ERC Advanced Grant Project MULTIMOD-26718.
J. Garnier and G. Papanicolaou   thank the Institut des Hautes \'Etudes
Scientifiques (IH\'ES) for its hospitality while this work was completed.

\end{document}